\documentclass[a4paper,10pt]{article}
\usepackage{epsfig}
\usepackage{color}
\usepackage{float}
\usepackage{hyperref}
\hypersetup{
colorlinks = true, 		
urlcolor    = cyan,  	
linkcolor  = blue,      
citecolor  = blue       
}

\usepackage{subcaption}
\usepackage{bm}
\usepackage{amsmath,amsthm,amssymb}
\usepackage{float}
\newtheorem{Teor}{Theorem}[section]

\newtheorem{Lem}{Lemma}[section]

\newtheorem{Rem}{Remark}[section]

\newcommand{\bea}{\begin{eqnarray}}
\newcommand{\eea}{\end{eqnarray}}
\newcommand{\beq}{\begin{equation}}
\newcommand{\eeq}{\end{equation}}
\newcommand{\beas}{\begin{eqnarray*}}
\newcommand{\eeas}{\end{eqnarray*}}

\newcommand{\dsp}{\displaystyle}

\newcommand{\A}{\boldsymbol{A}}

\title{Asymptotic behaviour of the semidiscrete FE approximations to weakly damped wave equations with minimal smoothness on initial data}
\author{ P. Danumjaya, Anil Kumar  and Amiya K. Pani\\ 
Department of Mathematics, BITS Pilani-KK Birla Goa Campus, \\ Zuarinagar, Goa-403726, India.\\ 
Email: danu@goa.bits-pilani.ac.in; anilpundir@goa.bits-pilani.ac.in and \\
amiyap@goa.bits-pilani.ac.in} 
\date{}
\allowdisplaybreaks
\begin{document}
\maketitle

\begin{abstract}
Exponential decay estimates of a general linear weakly damped wave equation are studied with decay rate lying in a range. Based on the $C^0$-conforming  finite element method to discretize spatial variables keeping temporal variable continuous, a semidiscrete system is analysed, and uniform decay estimates are derived with precisely the same decay rate as in the continuous case. Optimal error estimates with minimal smoothness assumptions on the initial data are established, which preserve exponential decay rate, and for a 2D problem, the maximum error bound is also proved. 
The present analysis is then generalized to  include the problems with non-homogeneous forcing function, space-dependent damping, and problems with compensator. It is observed that decay rates are improved with large viscous damping and compensator. Finally, some numerical experiments are performed to confirm our  theoretical findings.
\end{abstract}
	
\section{Introduction}
This paper deals with uniform exponential decay rates  for  the semidiscrete finite element solution of the following weakly damped wave equation:
\bea 
 u_{tt}(x,t) + \alpha \, u_{t}(x,t) -\nabla\cdot\big(\A(x) \nabla u(x)\big) + a_0(x) u(x)= 0, \; x \in \Omega, \; t > 0, \label{E1.1}
\eea
with initial conditions
\bea
u(x, 0) = u_0(x), \quad u_t(x, 0) = u_1(x), \; x \in \Omega, \label{E1.2}
\eea
and the boundary condition
\bea
u(x,t) = 0, \quad (x, t) \in \partial \Omega \times (0, \infty), \label{E1.3}
\eea
where $u_t = \frac{\partial u}{\partial t}$, $\Omega$ is a convex polygonal or polyhedral domain in $\mathbb{R}^d$ with boundary $\partial \Omega$,  and $\alpha$ is a fixed positive constant. Here, the matrix  $\A (x)$ is  a real symmetric and uniformly positive definite matrix for all $x\in \Omega$ and $a_0(x) \geq 0$ with coefficients  $\A$ and $a_0$ being smooth functions.


The equation (\ref{E1.1}) is known as the damped wave or telegraphers equation \cite{LCE, HET}, which arises in many applications such as acoustics, linear elasticity and  electro-magnetics, etc. Due to many applications, the damped wave equation has attracted significant interest in the literature. For the existence of a weak solution with regularity results using the Bubnov-Galerkin method  and weak compactness arguments, we may refer to \cite[Theorems 4.1-4.2 of Chapter II]{Temam} and \cite[Section 1.8]{BV}.

We use the standard notation for Sobolev spaces and their norms. In particular, let $L^2(\Omega)$ denote the space of square integrable functions on $\Omega$ with natural inner product $(\cdot, \cdot)$ and induced norm $\|\cdot\|$. For a nonnegative integer $k$, let $H^k$ denote the Hilbert Sobolev space $H^k(\Omega)$ with norm $\|\cdot\|_k$. Let
$$
H_0^1(\Omega) = \{ \phi \in H^1(\Omega) : \phi = 0 \;\;\mbox{on}\; 
\partial \Omega \}. 
$$ 
Denoting $A:H^2(\Omega) \cap H^1_0(\Omega) \longrightarrow L^2(\Omega)$ as a linear self-adjoint uniformly elliptic operator defined by
$$A \phi(x) = -\nabla\cdot\big(\A(x) \nabla \phi(x)\big) + a_0(x) \phi(x),$$
the problem \eqref{E1.1}-\eqref{E1.3} in abstract form is to 
seek  $u(t)\in D(A)$ for $t >0$ satisfying 
\begin{equation}\label{form-abstract}
u''(t) + \alpha u'(t) + Au(t)=0,\;t> 0,\;\;u(0)=u_0,\;\;\;u'(0)=u_1,
\end{equation}
where $u'(t)=\frac{du}{dt}$ and $u''(t)=\frac{d^2u}{dt^2}.$ With $A$ in \eqref{form-abstract} as a linear self-adjoint and positive definite operator  on $L^2(\Omega)$ with dense domain $D(A)=H^2(\Omega)\cap H^1_0(\Omega),$ we   define   $D(A^{r/2})$ as $\dot{H}^r=\dot{H}^r(\Omega),$ which is    a subspace of $H^r$ with norm $|v|_{r} =\|A^{r/2}\|.$ 
Essentially, for $r\geq 0$ and $r/2-1/4$ is not an integer, the space $\dot{H}^r=\{v\in H^r: A^{j} v= 0 \;\mbox{on}\;\partial\Omega, \; \mbox{for}\;j < r/2 -1/4\}$ and its norm is equivalent to $H^r$-norm.  Now, $H^1_0=D(A^{1/2})=\dot{H}^1_0$ and $\dot{H}^2= H^2\cap H^1_0.$

In the literature, explicit nonuniform decay rates  have been established  using a control-theoretic method for weakly damped 
 linear systems in Hilbert space, see, \cite{Russell-1}. The decay rate for the problem \eqref{E1.1} is given in terms of 
 the principal eigenvalue  of the operator ${A}$ and  the weak damping coefficient $\alpha>0$ in \cite[Proposition 1.2 
 of Chapter 4]{Temam}. In this paper, we have proved a better decay rate not only for the first energy, but also for higher 
 order energy. When the damped coefficient $\alpha=\alpha(x)>0,$ the decay rate involves the minimum and maximum of this
  coefficient and the  principal  eigenvalue of ${A}$ in \cite{Rau}. For related papers, see, also \cite{chen-1}-\cite{chen-2} and
   references, therein. In all these papers mentioned above, a large damping coefficient does not necessarily give rise to a large
   decay rate as it also depends on the principal  eigenvalue of the associated elliptic eigenvalue problem. Subsequently, Chen
   \cite{chen-3}  has developed and analysed improved decay rates by a new stabilization scheme that combines viscous damping
    and compensation. We shall discuss it in section 4.3 under generalizations.

In general, uniform decay property of  the continuous problem \eqref{E1.1} may not be preserved by  the approximate  solution when standard numerical schemes are applied. This may be due to existence of high frequency modes which  are only  weakly damped. Therefore, several stabilized methods have been developed and analysed, which give rise to uniform decay property of the semidiscrete-in-space schemes, keeping time variable continuous, see,  \cite{RTT} and \cite{Z-2005}  and references, therein. It is to be noted that mixed finite element methods are also employed to preserve uniform exponential decay property, see, \cite{HET}. This paper follows a different strategy to discuss uniform decay property of the semidiscrete scheme, when $C^0$-conforming finite element method is applied in the spatial direction. The key to the success of the present scheme is based on the energy arguments with the bound on the Poincar\'{e} type inequality, which provides a decay estimate in a range and similar to the decay rate as predicted by the continuous problem. However, the decay rate given by the present analysis may not be optimal and this is due  to non-conservative bounds in our estimates.
The main contributions of this paper are as follows. 
\begin{itemize}
\item [(i)] The first part of this paper focusses on the problem \eqref{E1.1}-\eqref{E1.3} and higher order in time regularity results are derived along with exponential decay properties using  energy arguments of  \cite[Theorems 4.1-4.2 of Chapter II]{Temam}  and  \cite[Section 1.8]{BV}. It is observed that the decay rate is calculated in a range involving the  damping parameter $\alpha$ and the first positive eigenvalue of the operator $A.$ In case, $u_0\in D(A^{(k)})$ and $u_1\in D(A^{(k-1/2)})$, the corresponding energy $2{\cal {E}}_{A^{(k)}}(u) =\|A^{(k-1/2)} u'\|^2 + \|A^{(k)} u\|^2$ decays exponentially with the same decay rate. In fact, it is observed that for large damping parameter, the decay rates may not be higher.  
\item [(ii)] Based on $C^0$-conforming finite element (FE) discretization in spatial variables keeping time variable continuous, a semidiscrete scheme is proposed, and uniform exponential decay estimates, which are uniform with respect to the discretizing parameter,  are derived. 
\item [(iii)] Optimal error estimates are established with minimal smoothness assumption on the initial data, that is, when $u_0\in H^3\cap H^1_0$ and $u_1\in H^2\cap H^1_0,$ which  have the same decay rate as observed for the semidiscrete solution. When $d=2$, the maximum norm estimate is also obtained with exactly  same decay rate.  
\item [(iv)] The analysis is then extended to include the nonhomogeneous problem and the problem with the space dependent viscous damping.
\item [(v)]  Decay rates are improved by the new stabilization method of combined with viscous damping and compensator for the semidiscrete solution.  Compared to (i)-(iii), the decay rates  can be made larger by choosing  large damping parameter and large compensator. 
\item [(vi)] Finally, several numerical experiments are conducted to confirm our theoretical findings.
\end{itemize}
Regarding (iii), earlier
 Rauch \cite{Rauch} has  initiated  optimal order of convergence, when $C^0$-conforming linear finite element  method is applied to approximate  a second order linear wave equation with minimal  smoothness on initial data, that is, $u_0\in H^3\cap H^1_0$ and $u_1=0,$  see also \cite{SP}, \cite{S}, \cite{KPY} and references, therein. Here, we emphasise that we have derived  results for the present problem with minimal smoothness on initial conditions $u_0$ and $u_1.$

An outline of this paper is as follows. In Section 2, we discuss weak formulation, regularity, and decay estimates for the continuous problem. Section 3 deals with the semidiscrete scheme. We establish decay estimates and optimal error estimates for the semidiscrete scheme. Section 4 is devoted to some generalizations involving  inhomogeneous problems, space dependent damping problems, and problems with  damping and compensator. Section 5 discuses  a  completely discrete scheme along with its energy conservation properties.  Finally,  several numerical experiments are conducted, whose results confirm  our theoretical findings.

\section{Weak formulation, Regularity results and Decay properties}
\setcounter{equation}{0}
This section deals with the weak formulation, some regularity results, and also the decay properties for the continuous problem.


We  now state   the following theorem  on  existence of a  unique weak solution, whose proof can be found in \cite[Theorem 1.1]{Lions-Strauss}, \cite[Theorem 4.1 of Chapter II]{Temam}.  
\begin{Teor}
Assume that $u_0 \in D(A)$ and $u_1 \in D(A^{1/2}).$ Then, the problem (\ref{E1.1})-(\ref{E1.3}) admits a unique strong solution $u$ satisfying 
\[ u \in L^{\infty}\left(0, T; D(A)\right), \quad  u' \in L^{\infty} \left(0, T; D(A^{1/2}) \right), \quad u'' \in L^{\infty} \left(0, T; L^2(\Omega) \right), 
\]
and 
\[ u'' + \alpha \, u' + A u = 0,\; \; \mbox{a.e.} \; t > 0 
\]
with 
\[ u(0) = u_0, \quad u'(0) = u_1.
\] 
\end{Teor}
Now, the bilinear form $a(\cdot,\cdot)$ on $V=D(A^{1/2})$ associated with $A$ is defined for $v,w\in D(A^{1/2})$ by
\begin{equation*}\label{bilinear-form}
a(v,w):=(A^{1/2} v, A^{1/2} w):= \left(\A(x) \nabla v,\nabla w \right) + \left( a_0 v, w\right).
\end{equation*}
Then, rewrite \eqref{form-abstract} as
\begin{eqnarray}
(u'',\chi) + \alpha (u',\chi) + a(u, \chi) = 0, \;\chi \in D(A^{1/2}), \label{E1.4} \\
u(0) = u_0, \quad \mbox{and} \quad u'(0) = u_1. \label{E1.5}
\end{eqnarray}

Let us recall here the following Poincar{\'e} inequalities for our subsequent use. For $v \in  D(A^{1/2})$ 
\begin{equation} \label{P1}
\|v \| \leq \frac{1}{\sqrt{\lambda_1}} \|A^{1/2} v\|, 
\end{equation}
and for $v \in  D(A) =H^2 \cap H_0^1$,
\begin{equation} \label{P2}
\|A^{1/2} v \| \leq \frac{1}{\sqrt{\lambda_1}} \| A v\|,
\end{equation}
where $\lambda_1$ is the principal eigenvalue of $A.$ 

\subsection{Decay Property}
This subsection focuses on the decay properties for the continuous problem (\ref{E1.1})-(\ref{E1.3}). Now, define the energy functional
\begin{eqnarray}
{\cal E}^{(1)}(u)(t) = \frac{1}{2} \left(\|u'\|^2 +\|A^{1/2} u\|^2\right). \label{E1.6}
\end{eqnarray}
\begin{Teor} \label{Th2.1}
For  $0 < \delta \leq \frac{1}{3} \min \left(\alpha, \frac{\lambda_1}{2 \alpha} \right),$ the solution $u$ of (\ref{E1.1})-(\ref{E1.3}) satisfies 
\begin{equation}
{\cal E}^{(1)}(u)(t) \leq 3 \, e^{-2 \delta t} \, {\cal E}^{(1)}(u)(0), \; t \geq 0. \end{equation}
\end{Teor}
\noindent
{\em Proof.} Set $\chi = u' + \epsilon u$ in (\ref{E1.4}) 
and then, rewrite the resulting equation as 
\begin{equation*}
	\frac{1}{2} \frac{d}{dt}\left(\|u'\|^2 + \|A^{1/2} u\|^2 + \epsilon \, \alpha \|u\|^2 \right) + \epsilon (u'', u) + \alpha \|u'\|^2 + \epsilon \|A^{1/2}  u\|^2 = 0. \label{E1.8}
\end{equation*}
A use of  the energy (\ref{E1.6})  with $\epsilon (u'',u) = \epsilon \frac{d}{dt}(u', u) - \epsilon \|u'\|^2$ 
shows
\begin{eqnarray}\label{e1.12}
\frac{d}{dt} \left( {\cal E}^{(1)}(u)(t)  +\epsilon (u', u) + \frac{\alpha \epsilon}{2}  \|u\|^2\right)  + \left((\alpha - \epsilon)  \|u'\|^2  + \epsilon \|A^{1/2}  u\|^2 \right) = 0.
\end{eqnarray}
Setting
$${\cal E}_{\epsilon}^{(1)}(u) (t)  = {\cal E}^{(1)}(u)(t) + \epsilon (u', u)  + \frac{\alpha \epsilon}{2}  \|u\|^2 \;\;
\mbox{and}\; \;F(t) = (\alpha - \epsilon) \|u'\|^2 + \epsilon \|A^{1/2}  u\|^2,
$$
rewrite (\ref{e1.12}) as
\begin{equation}\label{e1.13}
\frac{d}{d t} {\cal E}_{\epsilon}^{(1)}(u) (t) + F(t) = 0. 
\end{equation}
Observe that
\begin{eqnarray*}
	{\cal E}_{\epsilon}^{(1)}(u)(t) 
	\geq {\cal E}^{(1)}(u)(t) - \epsilon \|u'\| \|u\| + \frac{\alpha \epsilon}{2}  \|u\|^2.
\end{eqnarray*}
Using the Young's inequality, we obtain
\begin{eqnarray*}
	{\cal E}_{\epsilon}^{(1)}(u)(t) \geq {\cal E}^{(1)}(u)(t) - \epsilon \left(\frac{1}{2\alpha} \|u'\|^2 + \frac{\alpha}{2} \|u\|^2\right) + \frac{\alpha \epsilon}{2}  \|u\|^2 = {\cal E}^{(1)}(u)(t) - \frac{\epsilon}{2\alpha} \|u'\|^2, \label{E1.16}
\end{eqnarray*}
and hence, 
\begin{eqnarray*}
{\cal E}_{\epsilon}^{(1)}(u)(t) \geq \left(1- \frac{\epsilon}{\alpha}\right) {\cal E}^{(1)}(u)(t). \label{E1.17}
\end{eqnarray*}
Choose $\epsilon>0$ so that $\frac{\epsilon}{\alpha} \leq \frac{1}{2}$, i.e., $0 < \epsilon \leq \frac{\alpha}{2}$ 
to arrive at
\begin{equation} \label{E1.18}
{\cal E}_{\epsilon}^{(1)}(u)(t) \geq \frac{1}{2} {\cal E}^{(1)}(u)(t).
\end{equation}
Again, recall the definition of ${\cal E}_{\epsilon}^{(1)}(t)$ and use the Cauchy-Schwarz and Young's inequalities to find that 
\begin{eqnarray*}
	{\cal E}_{\epsilon}^{(1)}(u)(t)  
	\leq \frac{3}{2}{\cal E}^{(1)}(u)(t)  +\left(\frac{\epsilon}{2\alpha} -\frac{1}{4}\right) \|u'\|^2 -\frac{1}{4} \|A^{1/2} u\|^2 +   \epsilon \alpha \|u\|^2, \label{E1.19}
\end{eqnarray*}
and a use of the Poincar{\'e} inequality (\ref{P1}) yields
\begin{eqnarray*} \label{E1.20}
	{\cal E}_{\epsilon}^{(1)}(u)(t) \leq \frac{3}{2}{\cal E}^{(1)}(u)(t)  +\left(\frac{\epsilon}{2\alpha} -\frac{1}{4}\right) \|u'\|^2 + \left(\frac{\epsilon \alpha}{\lambda_1}-\frac{1}{4}\right) \|A^{1/2} u\|^2.	
\end{eqnarray*} 
In order to derive an estimate of the form ${\cal E}_{\epsilon}^{(1)}(u)(t) \leq  \frac{3}{2}{\cal E}^{(1)}(u)(t)$, we must have 
\begin{equation}
\frac{\epsilon}{2\alpha} -\frac{1}{4} \leq 0 \quad \mbox{and} \quad \frac{\epsilon \alpha}{\lambda_1}-\frac{1}{4} \leq 0, \label{E1.21}
\end{equation}
that is, set 
$$ 
0 < \epsilon \leq \frac{1}{2} \min \left(\alpha, \frac{\lambda_1}{2 \alpha}\right).
$$
Therefore, combining (\ref{E1.18}) and (\ref{E1.21}), we obtain
\begin{equation}\label{E1.22}
\frac{1}{2}{\cal E}^{(1)}(u)(t) \leq {\cal E}_{\epsilon}^{(1)}(u)(t) \leq  \frac{3}{2}{\cal E}^{(1)}(u)(t),
\end{equation}
provided $0 < \epsilon \leq \frac{1}{2} \min \left(\alpha, \frac{\lambda_1}{2 \alpha}\right)$. A use of the definition of $F$ and ${\cal E}^{(1)}(u)(t)$ shows
\begin{equation*}
F(t) = (\alpha - 2 \epsilon) \|u'\|^2 +2\epsilon {\cal E}^{(1)}(u)(t).
\end{equation*}
Note that for $0 < \epsilon = \min\left(\frac{\alpha}{2}, \frac{\lambda_1}{4 \alpha}\right)$, $\alpha -2\epsilon \geq 0$, there holds
\begin{equation*}\label{estimate:F}
F(t) \geq 2\epsilon {\cal E}^{(1)}(u)(t) \geq  \frac{4 \epsilon}{3}  {\cal E}_{\epsilon}^{(1)}(u)(t).
\end{equation*}
Thus, from (\ref{e1.12}), it follows that
\begin{equation*} \label{E1.22-1}
\frac{d}{dt} {\cal E}_{\epsilon}^{(1)}(u)(t) + \frac{4 \epsilon}{3} {\cal E}_{\epsilon}^{(1)}(u)(t) \leq \frac{d}{dt} {\cal E}_{\epsilon}^{(1)}(u)(t) + F(t) = 0,
\end{equation*}
and then, an integration with respect to $t$ shows
\begin{eqnarray*}
	{\cal E}_{\epsilon}^{(1)}(u)(t) \leq e^{- \frac{4}{3}\epsilon t} \,{\cal E}_{\epsilon}^{(1)}(u)(0) \leq  \frac{3}{2} \, e^{- \frac{4}{3}\epsilon t}\, {\cal E}^{(1)}(u)(0).
\end{eqnarray*}
A use of (\ref{E1.22}) yields
\begin{equation*}\label{E1.23}
{\cal E}^{(1)}(u)(t) \leq 3 \, e^{-2\delta t} \, {\cal E}^{(1)}(u)(0),
\end{equation*}
where $\delta \in \left(0, \frac{1}{3} \min \left(\alpha, \frac{\lambda_1}{2 \alpha} \right) \right)$. This completes the rest of the proof. \hfill$\blacksquare$ \\

The next theorem is on higher order time derivatives of energy.
\begin{Teor} \label{Th2.2}
For $0 < \delta \leq \frac{1}{3} \min \left(\alpha, \frac{\lambda_1}{2 \alpha} \right),$ the solution $u$ of (\ref{E1.1})-(\ref{E1.3}) satisfies 

\begin{equation*}
{\cal E}^{(j)}(u)(t) \leq 3 \,e^{-2 \delta t} \, {\cal E}^{(j)}(u)(0), \; j =  2, 3,\ldots, \; t \geq 0, \label{Neq2}
\end{equation*}
where 
$$
{\cal E}^{(j)}(u)(t) = \frac{1}{2} \left(\|u^{(j)}\|^2 + \| A^{1/2}  u^{(j-1)}\|^2 \right),
$$
and $u^{(j)}$ stands for $j^{\mbox{th}}$ time derivative of $u$.
\end{Teor}
\noindent
{\em Proof.} 
On differentiating $j\geq 2$ times the equation \eqref{E1.4}, we easily obtain
\begin{equation}\label{E1.24}
(u^{(j+1)},\chi) + \alpha (u^{(j)}, \chi) + a ( u^{(j-1)}, \chi) = 0,\;\;\;\chi \in D(A^{1/2}).
\end{equation}
With $w = u^{(j-1)}$, $w^{(0)} =w,$  $w^{(1)}= w'$ and $w^{(2)}=w''$ in (\ref{E1.24}) , we arrive at an equation  \eqref{E1.4} now in $w$ as 
\begin{equation*}\label{E1.24-w}
(w'',\chi) + \alpha (w', \chi) + a (w ,  \chi) = 0, \;\;\;\chi \in D(A^{1/2}).
\end{equation*}
Therefore, we repeat the argument as in Theorem ~\ref{Th2.1} and obtain result in terms of $w.$ Then, writing in terms of $u$, we complete the rest of the proof. \hfill $\blacksquare$

\begin{Teor}
For $0 < \delta \leq \frac{1}{3} \min \left(\alpha, \frac{\lambda_1}{2 \alpha} \right),$ the solution $u$ of (\ref{E1.1})-(\ref{E1.3}) satisfies 
\begin{equation*}
{\cal E}_A^{(1)}(u)(t) \leq 3 \,e^{-2 \delta t} \, {\cal E}_A^{(1)}(u)(0), \; j = 1, 2, 3, \ldots, \; t \geq 0, \label{Neq1}
\end{equation*}
where
\begin{equation*}
{\cal E}_A^{(1)}(u)(t) = \frac{1}{2} \left(\|A^{1/2} u'(t)\|^2 + \|A u(t)\|^2 \right).
\end{equation*}
\end{Teor}
\noindent
{\em Proof.} The analysis closely follows the proof technique of Theorem \ref{Th2.1}. Forming an inner product between (\ref{E1.1}) and $A (u' + \epsilon u)$,  then rewrite it to arrive at
\begin{equation*} 
\frac{d}{d t} {\cal E}_{A, \epsilon}^{(1)}(u)(t) + (\alpha - \epsilon) \|A^{1/2} u'(t)\|^2  + \epsilon \|A u(t)\|^2 = 0,
\end{equation*}
where
$$
{\cal E}_{A, \epsilon}^{(1)}(u)(t) = {\cal E}_A^{(1)}(u)(t) + \epsilon \left(A^{1/2} u'(t), A^{1/2} u(t)\right) + \frac{\alpha\, \epsilon}{2} \|A^{1/2} u(t)\|^2.
$$  
We now proceed exactly in the proof technique of Theorem \ref{Th2.1} by replacing ${\cal E}^{(1)}$ by ${\cal E}_A^{(1)}, \; {\cal E}_{\epsilon}^{(1)}$ by ${\cal E}_{A, \epsilon}^{(1)}$ and using Poincar{\'e} inequality (\ref{P2}) to arrive at
$$
{\cal E}_A^{(1)}(u)(t) \leq 3 \, e^{-2\delta t} \, {\cal E}_A^{(1)}(u)(0),
$$
whenever $0 < \delta \leq \frac{1}{3} \min \left(\alpha, \frac{\lambda_1}{2 \alpha} \right)$. This completes the rest of the proof. \hfill $\blacksquare$

\begin{Rem}
Since 
$$
\|Au(t)\|^2 \leq 3 e^{-2 \delta t} {\cal E}^{(1)}_A (u) (0) \leq \frac{3}{2} \left(\|A^{1/2} u_1 \|^2 + \| A u_0\|^2 \right) \leq \frac{3}{2} 
\left(\| u_1\|_{1}^2 + \| u_0\|^2_2 \right).
$$
A use of elliptic regularity yields $\|\Delta u(t)\| \geq C_R \|u(t)\|_2$ with the Sobolev embedding result   shows
\begin{equation*} \label{New2.29}
\|u(t)\|_{L^{\infty}} \leq C \|u(t)\|_2 \leq C e^{-\delta t} \left(\|u_1\|_1 + \|u_0\|_2 \right).
\end{equation*}
\end{Rem}

\begin{Rem}
Following the proof technique of Theorem 2.4, the following result 
\begin{equation*} \label{New2.30}
{\cal E}_A^{(j)}(u)(t) \leq 3 \,e^{-2\delta t} \, {\cal E}_A^{(j)}(u)(0),
\end{equation*}
where 
\begin{equation*} \label{New2.31}
{\cal E}_A^{(j)}(u)(t) = \frac{1}{2} \left(\|A^{1/2} u^{(j)}(t)\|^2 + \|A u^{(j-1)}(t)\|^2 \right),
\end{equation*}
can be proved by using induction hypothesis. \\

Assume that  $u_0\in D(A^{(k)})$ and $u_1\in D(A^{(k-1/2)})$ for $k>1.$  Then, following the arguments in Theorem 2.2-2.4 and using induction, there holds:
$${\cal E}_{A^{(k)}}^{(j)} (u)(t) \leq 3 \,e^{-2\delta t} \, {\cal E}_{A^{(k)}}^{(j)}(u)(0),$$
where ${\cal E}_{A^{(k)}}^{(j)} (u)(t) =\frac{1}{2} \left(\|A^{(k-1/2)} u^{(j)}(t)\|^2 + \|A^{(k)} u^{(j-1)}(t)\|^2 \right).$
\end{Rem}

\section{Semidiscrete scheme}
\setcounter{equation}{0}
This section analyses the semidiscrete method for the problem (\ref{E1.1})-(\ref{E1.3})  
and discusses the decay rates along with the optimal error estimates. 

Let $\{S_h^0\}_{h>0}$ be a  family of subspaces  of $H^1_0$ with the following approximation property:
\begin{equation}\label{approx-property}
\inf_{\chi \in S_h^0} \left( \|v-\chi\| + h \|v-\chi\|_1\right) \leq h^r\;\|v\|_r,\;\;\mbox{for}\;v\in H^r\cap H^1_0.
\end{equation}
The semidiscrete formulation is to find  $u_h: [0,\infty) \rightarrow S_h^0$ such that
\begin{eqnarray}
(u_h''(t),\chi) + \alpha (u_h'(t),\chi) + a (u_h(t), \chi) = 0, \;\chi \in S_h^0, \label{E2.1} \\
u_h(0) = u_{0,h}, \; \mbox{and} \; u_h'(0) = u_{1,h}, \label{E2.2}
\end{eqnarray}
where $u_{0,h}$ and $u_{1,h}$ are appropriate approximations of $u_0$ and $u_1,$ respectively, in $S_h^0$
to be defined later. Since $S_h^0$ is finite dimensional, (\ref{E2.1}) gives rise to a system of linear ODEs. An application of the Picard's theorem yields the existence of a unique discrete solution $u_h(t)\in S^0_h,$ for all $ t \in (0, \infty)$.

Let us first define a discrete counterpart  $A_h : S_h^0 \mapsto S_h^0$   of the operator $A$  as 
\begin{equation}\label{discrete-A}
(A_h v_h, \chi) = a (v_h,\chi)\;\;\forall v_h, \chi \in S^0_h.
\end{equation}
Then, we rewrite \eqref{E2.1} as 
\begin{equation}\label{discrete-problem-2}
u_h'' + \alpha u_h'+  A_h u_h =0,\;\; t>0.
\end{equation}

\subsection{Decay Property}
This subsection discusses the decay estimates for the solution of semidiscrete equation.
Now, define the energy functional as
\begin{eqnarray*}
{\cal E}_h^{(1)}(u_h)(t) = \frac{1}{2} \left(\|u_h'(t)\|^2 +\| A_h^{1/2} u_h(t)\|^2\right), 
\end{eqnarray*}
where $\|A_h^{1/2}u_h\|^2:= a( u_h,u_h).$
\begin{Teor}
For $0 < \delta \leq \frac{1}{3} \min \left(\alpha, \frac{\lambda_1}{2 \alpha} \right)$, the solution $u_h$ of (\ref{E2.1})-(\ref{E2.2}) satisfies the following decay property
\begin{equation*}
{\cal E}_h^{(1)}(u_h)(t) \leq 3 \, e^{-2 \delta t} \,{\cal E}_h^{(1)}(u_h)(0),  \; t \geq 0.
\end{equation*}
\end{Teor}
\noindent
{\em Proof.} A use of $\chi = u_h'(t) + \epsilon \, u_h(t)$ in (\ref{E2.1}) yields
\begin{eqnarray*}
\frac{d}{dt} \left( {\cal E}_h^{(1)}(u_h)(t) + \epsilon (u_h', u_h) + \frac{\alpha \epsilon}{2} \|u_h(t)\|^2 \right) + (\alpha - \epsilon) \|u_h'(t)\|^2 + \epsilon \|A_h^{1/2} u_h(t)\|^2 = 0.
\end{eqnarray*}
Since $u_h(t) \in S_h^0 \subset H_0^1,$  then by Poincar{\'e} inequality \eqref{P1}
\begin{equation*} \label{Ph1}
\|u_h(t)\| \leq \frac{1}{\sqrt{\lambda_1}} \|A^{1/2}_h u_h(t)\|.
\end{equation*}
We then proceed exactly like the proof of the Theorem \ref{Th2.1} replacing $u$ by $u_h$ to obtain
\begin{equation*}
{\cal E}_h^{(1)}(u_h)(t) \leq 3 \, e^{-2 \delta t} \, {\cal E}_h^{(1)}(u_h)(0), \; t \geq 0.
\end{equation*}
This completes the rest of the proof. \hfill $\blacksquare$

\begin{Teor}
For $0 < \delta \leq \frac{1}{3} \min \left(\alpha, \frac{\lambda_1}{2 \alpha} \right)$, the solution $u_h$ of (\ref{E2.1})-(\ref{E2.2}) satisfies
\begin{equation*}
{\cal E}_h^{(j)}(u_h)(t) \leq 3 \,e^{-2 \delta t} \, {\cal E}_h^{(j)}(u_h)(0), \; j = 2, 3,\ldots, \; t \geq 0, 
\end{equation*}
where 
\begin{eqnarray*}
{\cal E}_h^{(j)}(u_h)(t) = \frac{1}{2} \left(\|u_h^{(j)}(t)\|^2 +\| A_h^{1/2} u_h^{(j-1)}(t)\|^2\right).
\end{eqnarray*}
\end{Teor}
\noindent
{\em Proof.} We prove the result ${\cal E}_h^{(j)}(u_h)$ by using the induction hypothesis. Assume that the result is true for $j - 1,$ that is,
$$
{\cal E}_h^{(j-1)}(u_h)(t) \leq 3 \, e^{-2 \delta t} \, {\cal E}_h^{(j-1)}(u_h)(0).
$$
We now consider
\begin{equation*}
(u_h^{(j+1)},\chi) + \alpha (u_h^{(j-1)}, \chi) +a  ( u_h^{(j-1)},  \chi) = 0.
\end{equation*}
Choose $w_h = u_h^{(j-1)}$ and $\chi = w_h' + \epsilon \, w_h$ and follow similar steps like proof of Theorem \ref{Th2.2} replacing ${\cal E}(u)$ by ${\cal E}_h(u_h)$ to obtain 
\begin{equation*}
{\cal E}_h^{(j)}(u_h)(t) \leq 3 \, e^{-2 \delta t}\,{\cal E}_h^{(j)}(u_h)(0), \; j = 2, 3,\ldots, \; t \geq 0. 
\end{equation*}
This completes the rest of the proof. \hfill $\blacksquare$
\begin{Rem}
If $\|A_h^{1/2} u_{0 h}\| \leq C \| u_0\|_{1}$ and $\|u_{1 h} \| \leq C \|u_1\|$ then 
$$
\|A^{1/2}_h u_h(t) \| \leq C e^{-\delta t} \left(\|u_{1 h}\| + \| A^{1/2}_h u_{0 h}\| \right) \leq C e^{-\delta t} \left(\|u_1\| + \| u_0\|_1 \right).
$$
Observe that using coercivity property of the bilinear form  $\|A_h^{1/2}u_h\|^2 = a( u_h,u_h) \geq \alpha_0 \|\nabla u_h\|^2,$ we arrive at
$$
\|\nabla u_h(t) \| \leq C e^{-\delta t} \left(\|u_1\| + \| u_0\|_1 \right).
$$
As a consequence of the Sobolev embedding for $d=2$, see, \cite{VT}, we obtain
$$
\|u_h(t)\|_{L^{\infty}} \leq C \left(\log \left(\frac{1}{h} \right)\right) \|\nabla u_h(t)\| \leq C \left(\log \left(\frac{1}{h} \right)\right) e^{-\delta t} \left(\|u_1\| + \|\nabla u_0\| \right).
$$
\end{Rem}

\begin{Teor}
For $0 < \delta \leq \frac{1}{3} \min \left(\alpha, \frac{\lambda_1}{2 \alpha} \right)$, and a positive constant $C,$ the solution $u_h$ of (\ref{E2.1})-(\ref{E2.2}) satisfies
\begin{equation*}
{\cal E}_{A_h}^{(1)}(u_h)(t) \leq 3 \,e^{-2 \delta t} \, {\cal E}_{A_h}^{(1)}(u_h)(0), \; t \geq 0,
\end{equation*}
where
$$
{\cal E}_{A_h}^{(1)}(u_h)(t) = \frac{1}{2} \left(\|A^{1/2}  u'_h(t)\|^2 + \|A_h u_h(t) \|^2 \right).
$$
\end{Teor}
\noindent
{\em Proof.} 
Forming inner product between equation \eqref{discrete-problem-2} and $A_h u'_h + \epsilon u_h$ to obtain
\begin{eqnarray*}
\frac{d}{d t} \left( {\cal E}_{A_h}^{(1)}(u_h) + \epsilon (A^{1/2}_h u'_h, A^{1/2}_hu_h) + \frac{\alpha \, \epsilon}{2} \|A_h^{1/2} u_h\|^2 \right) + (\alpha - \epsilon) \|A_h^{1/2}  u'_h\|^2 + \epsilon \|A_h u_h\|^2 = 0.
\end{eqnarray*}   
We then proceed in a similar manner exactly like the proof of Theorem \ref{Th2.1} and using for $v_h \in S_h^0$, Poincar{\'e} inequality (\ref{Ph1}), 
\begin{eqnarray*}
\|A_h^{1/2}  v_h\|^2 = (A_h v_h, v_h) &\leq& \|A_h v_h\| \, \|v_h\| \\
&\leq& \frac{1}{\sqrt{\lambda_1}} \|A_h v_h\| \, \|A_h^{1/2} v_h\|,
\end{eqnarray*}
that is,
$
\|A_h^{1/2}  v_h\| \leq \frac{1}{\sqrt{\lambda_1}} \|A_h v_h\|
$
and obtain
\begin{equation*}
{\cal E}_{A_h}^{(1)}(u_h)(t) \leq C \, e^{-2 \delta t} \, {\cal E}_{A_h}^{(1)}(u_h)(0), \; t \geq 0.
\end{equation*}
This completes the rest of the proof. \hfill $\blacksquare$

\subsection{Error estimates.} 
This subsection deals with optimal error estimates for the
semidiscrete scheme.  Throughout this subsection, we shall use $r=2$, that is, $S_h^0$ consisting of  $C^0$-conforming piecewise linear elements and for general $r>2$, all the ensuing results hold under assumptions of higher regularity on the exact solution. 

Let $R_h u$ be the elliptic projection of $u$ defined by
\begin{eqnarray}
a(u - R_h u, \chi) = 0, \; \forall \; \chi \in S_h^0. \label{E2.4}
\end{eqnarray}  
We split the error as 
\begin{equation*} \label{EqnErr}
e := u - u_h = (u - R_h u) + (R_h u - u_h) := \eta + \theta.
\end{equation*}     
Note that $a (\cdot,\cdot)$ satisfies the boundedness and coercivity properties.
Setting $\eta = u - R_h u$, the following estimates are easy
to obtain
\begin{eqnarray}
\|\eta\|_j + \|\eta_t\|_j &\leq& C h^{r+1-j} \left(\sum_{m=0}^{1}
					      \left\|\frac{\partial^m
					      u}{\partial t^m} \right\|_{r+1}
					     \right), \; j=0,1. \label{E2.5} 
\end{eqnarray}
For details, see, \cite{Bre}.

We subtract the equation (\ref{E2.1}) from (\ref{E1.4}), and using the elliptic projection (\ref{E2.4}), we obtain the error equation in $\theta$ as
\begin{eqnarray}
(\theta'', \chi) + \alpha (\theta', \chi) +a (\theta, \chi) = -(\eta'', \chi) - \alpha (\eta', \chi), \forall
\, \chi \in S_h^0. \label{E2.6}  
\end{eqnarray}
\begin{Lem}\label{Lem:3.1}
Let $\theta$ satisfy (\ref{E2.6}). 
Then, there holds  for small $\delta_0>0$
\begin{eqnarray*}
{\cal E}_h^{(1)} (\theta)(t) \leq 3 e^{-2 \delta t (1-\delta_0)} {\cal E}_h^{(1)}(\theta)(0) + \frac{1}{\delta_0} \left( \frac{2}{\alpha} + \frac{\alpha}{2\lambda_1} \right)  \int_0^t e^{-2 \delta  (1-\delta_0) (t - s)} \left(\|\eta''\|^2 + \|\eta'\|^2 \right) \, ds,
\end{eqnarray*}
where ${\cal E}_h^{(1)} (\theta)(t):=\frac{1}{2} \big(\|\theta'(t)\|^2 + \|A_h^{1/2} \theta(t)\|^2\big).$
\end{Lem}
\noindent
{\bf Proof.} 
Choosing $\chi = \theta' + \epsilon \theta$ in (\ref{E2.6}), we note that $\epsilon (\theta'', \theta) = \epsilon \frac{d}{d t} (\theta', \theta) - \epsilon \|\theta'\|^2$ and then, setting 
\begin{equation*}
{\cal E}^{(1)}_{h,\epsilon}(\theta) = {\cal E}_h^{(1)}(\theta) + \epsilon (\theta', \theta) + \frac{1}{2} \alpha \, \epsilon \|\theta(t)\|^2,
\end{equation*} 
and 
\begin{eqnarray*}
F(t) &=& \left(\alpha - \epsilon \right) \|\theta'(t)\|^2 + \epsilon \|A^{1/2} \theta\|^2\\
&\;\;\;\;=&  2 \epsilon {\cal E}_h^{(1)}(\theta) (t) + \left( \alpha - 2 \epsilon \right) \|\theta'\|^2,
\end{eqnarray*}
we now arrive applying the Cauchy-Schwarz inequality with the Young's inequality, Poinc\'{a}re inequality \eqref{Ph1}  and  for some $\delta_0>0$ at 
\begin{eqnarray*}
\frac{d}{d t} {\cal E}^{(1)}_{h,\epsilon} (\theta)(t) + F(t) &=& -(\eta'' + \alpha \eta', \theta') -\epsilon(\eta'' + \alpha \eta', \theta)\nonumber\\
&\leq& \left(\|\eta''\| +\| \eta\|\right)\;\|\theta'\| + \frac{\epsilon}{\sqrt{\lambda_1}}\;\left(\|\eta''\| +\| \eta\|\right)\;\|A_h^{1/2}\theta\| \nonumber\\
&\leq& \frac{1}{2\delta_0}\left( \frac{1}{\alpha-\epsilon} + \frac{\epsilon}{\lambda_1}\right) \;(\|\eta''\|^2 +\| \eta\|^2)+ \delta_0 \;F(t).
 \label{EQ3.18}
\end{eqnarray*}
With $
0 < \delta \leq \frac{1}{2} \min \left( \alpha, \frac{\lambda_1}{2 \alpha} \right)
$, it follows that $(\alpha - 2 \epsilon) \geq 0$ and 
$$
F(t) \geq 2 \epsilon \,{\cal E}^{(1)}_{h}(\theta)(t) \geq \frac{4}{3} \epsilon \, {\cal E}_{h,\epsilon}^1 (\theta)(t).
$$
On substitution, we arrive with $\epsilon \leq \alpha/2$ at
\begin{equation}
\frac{d}{d t} {\cal E}^{(1)}_{h,\epsilon} (\theta)(t)  + \frac{4}{3} \epsilon  (1-\delta_0) \;{\cal E}_{h,\epsilon}^1 (\theta)(t) \leq \frac{1}{2\delta_0} \left( \frac{2}{\alpha} + \frac{\alpha}{2\lambda_1} \right) \;\left(\|\eta''(t)\|^2 +  \|\eta'(t)\|^2\right). \label{EQ3.19}
\end{equation} 
We rewrite the equation (\ref{EQ3.19}) as
\begin{eqnarray*}
\frac{d}{dt} \left(e^{\frac{4}{3} \epsilon (1-\delta_0) t} \,{\cal E}^{(1)}_{h,\epsilon} (\theta)(t) \right) \leq 
\frac{1}{2\delta_0} \left( \frac{2}{\alpha} + \frac{\alpha}{2\lambda_1} \right) \;e^{\frac{4}{3} \epsilon (1-\delta_0) t} \left(\|\eta''(t)\|^2 +  \|\eta'(t)\|^2\right).
 \label{EQ3.20}
\end{eqnarray*} 
On integration from $0$ to $t$, it follows that
\begin{eqnarray*}
{\cal E}^{(1)}_{h,\epsilon} (\theta)(t) &\leq& e^{-\frac{4}{3} \epsilon (1-\delta_0) t} {\cal E}^{(1)}_{h,\epsilon} (\theta)(0) 
+\frac{1}{2\delta_0} \left( \frac{2}{\alpha} + \frac{\alpha}{2\lambda_1} \right) \int_0^t e^{-\frac{4}{3} \epsilon (1-\delta_0)(t-s)} \big(\|\eta''(s)\|^2 + \|\eta(s)\|^2\big) \, ds.
\end{eqnarray*}
With $2 \delta = \frac{4}{3} \epsilon (1-\delta_0), $ that is, $ \delta = \frac{2}{3} \epsilon (1-\delta_0)$ and using ${\cal E}^{(1)}_{h \epsilon}$ in terms of ${\cal E}^{(1)}_{h}$, we complete the rest of the proof. \hfill $\blacksquare$ 

\begin{Rem} 
When $u_{0h} = R_h u_0,$ then $\theta(0) = 0$ and therefore, 
$$
{\cal E}_h^{(1)}(\theta)(0) = \frac{1}{2} \|\theta'(0)\|^2.
$$
With $u_{1 h}$ either $L^2$-projection or interpolant of $u_1$ in $S_h^0$, we obtain
$$
{\cal E}_h^{(1)}(\theta)(0) \leq C h^4 \|u_1\|_2^2.
$$
Therefore, we arrive at the following superconvergent result for $\|A^{1/2}_h \theta(t)\|$
\begin{eqnarray*}
\|\theta'(t)\|^2 + \|A_h^{1/2} \theta(t)\|^2 \leq  C h^4 e^{-2 \delta (1-\delta_0) t} \left(
 \|u_1\|^2_2 + \int_0^t e^{2 \delta (1-\delta_0) s} \left(\|u''(s)\|_2^2 + \|u'(s)\|_2^2 \right) \, ds \right).
\end{eqnarray*}
Since from  Remark 2.3 with $k=1$ and $j=3,$ there holds 
\begin{eqnarray*}
\|u''(t)\|_2^2 &\leq& 6 \,e^{-2\delta t} \,{\cal E}_{A^{(1)}}^{(3)}(0)\nonumber\\
                     &\leq& 3\,e^{-2\delta t} \,\Big( \|A^{1/2}u^{(3)}(0)\|^2 + \|A u^{(2)}(0)\|^2\Big),
\end{eqnarray*}
and
\begin{eqnarray*}
\|u'(t)\|_2^2 &\leq& 6 \,e^{-2\delta t} \,{\cal E}_{A^{(1)}}^{(2)}(0)\nonumber\\
                     &\leq& 3\,e^{-2\delta t} \,\Big( \|A^{1/2}u''(0)\|^2 + \|A u_1\|^2\Big).
\end{eqnarray*}
A use of  $u'' (0) =-\alpha u_1 - Au_0$  with $u^{(3)}(0)= -\alpha u''(0)-A u_1= (\alpha^2-A) u_1 + A u_0$ implies
\begin{eqnarray*}
\|u''(t)\|_2^2  &\leq& C\,e^{-2\delta t} \,\Big( \|A^{3/2} u_1 \|^2+\|A^{2} u_0\|^2\Big) \nonumber\\
&\leq& C\;e^{-2\delta t}  \Big(\|u_0\|^2_{4} + \|u_1\|_{3}^2 \Big),
\end{eqnarray*}
and
\begin{eqnarray*} \label{estimate:u-t}
\|u'(t)\|_2^2  &\leq& C\,e^{-2\delta t} \,\Big( \|A u_1 \|^2+\|A^{3/2} u_0\|^2\Big) \nonumber\\
&\leq& C\;e^{-2\delta t}  \Big(\|u_0\|^2_{3} + \|u_1\|_{2}^2 \Big).
\end{eqnarray*}
Hence, we obtain the  following superconvergence result
\begin{equation}\label{supercgt}
\|\theta'(t)\|^2 + \|A^{1/2}_h \theta(t)\|^2 \leq C h^4 \,(1+t)\; e^{-2 \delta (1-\delta_0) t} \Big( \|u_0\|^2_{4} + \|u_1\|_{3}^2 \Big).
\end{equation}
As a by-product and using triangle inequality with \eqref{E2.5}, there hods
\begin{equation}\label{estimate-e-t}
\|u'(t) - u_h'(t) \|^2 \leq C  (1+t)\; h^4 e^{-2 \delta (1-\delta_0) t} \;\Big( \|u_0\|^2_{4} + \|u_1\|_{3}^2 \Big).
\end{equation}
\end{Rem}
\noindent
When $u_{0 h}$ is chosen as $L^2$-projection or an interpolant, then 
$$
{\cal E}^{(1)}(\theta)(0) \leq C h^2 \left( \|u_0\|_2^2 + \|u_1\|_1^2 \right),
$$
and hence, using the coercivity of the bilinear form  $a(\cdot,\cdot),$ 
 we find that
$$
\alpha_0 \| \nabla\theta(t)\|^2 \leq \|A_h^{1/2} \theta\|^2 \leq C h^2 \, e^{-2 \delta (1-\delta_0) t}\;\left( \|u_0\|_2^2 + \|u_1\|_1^2 + \int_0^t \left(\|u''\|_1^2 + \|u'\|_2^2 \right) \right) \, ds. \\
$$ 
A use of 
$$\|u''(s) \|_1^2 \leq C \|A^{1/2} u''(s)\|^2\leq 6\;e^{-2\delta s}  {\cal E}_{A^{(1)}}^{(2)}(0) \leq C\;e^{-2\delta s}\;
\Big(\|u_0\|^2_{3} + \|u_1\|_{2}^2 \Big),
$$
shows the following optimality error estimate.
\begin{Teor} With either 
$u_{0 h}$ and $u_{1 h},$ respectively, as interpolant  or $L^2$ projections of of $u_0$ and $u_1,$  there holds the following optimal error estimate for small $\delta_0>0$
\begin{eqnarray*}
\|\nabla (u - u_h)(t) \|^2
\leq  C h^2 (1+t) \;e^{-2 \delta (1-\delta_0) t} \Big(\|u_0\|^2_{3} + \|u_1\|_{2}^2 \Big).
\end{eqnarray*}
\end{Teor}
As a consequence of superconvergent result of $\|\nabla \theta(t)\|$ in \eqref{supercgt}, we apply the Sobolev embedding lemma for $d=2$, (see, \cite{VT})  to obtain 
\begin{eqnarray*}
\|\theta(t)\|_{L^{\infty}} \leq C \left( \log \left(\frac{1}{h} \right) \right) \|\nabla \theta(t)\| \leq C \left( \log \left(\frac{1}{h} \right) \right) (1+t)^{1/2} e^{- \delta (1-\delta_0) t} h^2 \big(\|u_0\|_{4} + \|u_1\|_3\big).
\end{eqnarray*}
Since 
$$
\|\eta(t)\|_{L^{\infty}} \leq C h^2 \left( \log \left(\frac{1}{h} \right) \right) \|u(t)\|_{W^{2, \infty}} \leq C h^2 \left( \log \left(\frac{1}{h} \right) \right) \big( \|u_0\|_{W^{2, \infty}} + \|u_1\|_{W^{1,\infty}}\big),
$$
then, for  $d=2$  and small $\delta_0>0$ there holds
\begin{eqnarray*}
\|u(t) - u_h(t)\|_{L^{\infty}} \leq C h^2 \left( \log \left(\frac{1}{h} \right) \right) (1+t)^{1/2}  e^{- \delta (1-\delta_0) t} \big(\|u_0\|_{4} + \|u_1\|_3\big), \label{Neweqn}
\end{eqnarray*}
provided $\|u(t)\|_{W^{2, \infty}} = O \left( e^{- \delta t} \right)$. 

From the superconvergence result for $\|\nabla \theta(t)\|$ in \eqref{supercgt}, one obtains estimate of $\|\theta(t)\|$, but with the assumption of higher regularity, that is, $u_0\in H^4\cap H^1_0$ and $u_1\in H^3\cap H^1_0$
and only with $u_{0 h} = R_h u_0$. 

Below, we directly deduce using a modified version of Baker's arguments \cite{B-1976}, an  optimal error estimate of $\|u(t) - u_h(t)\|$ 
and $u_{0 h}$ as $L^2$-projection or interpolant of $u_0$ onto $S_h^0$. 

\begin{Teor}
Let $u$ and $u_h$ be a solution of (\ref{E1.4}) and (\ref{E2.1}), respectively. Then, there exists a positive constant $C$ independent of $h$ such that 
$$
\|u(t) - u_h(t)\| \leq C h^2 \, (1+t)^{1/2}\, e^{-\delta (1-\delta_0) t} \big(\|u_0\|_{3} + \|u_1\|_2\big).
$$
\end{Teor}

\noindent
{\bf Proof.} Integrate (\ref{E2.6}) with respect to $t$ and obtain
\begin{eqnarray} \label{En3.26}
(\theta'(t), \chi) + \alpha \, (\theta(t), \chi) +a (\hat{\theta}(t),  \chi) = (e'(0), \chi) + \alpha (e(0), \chi) - (\eta', \chi) - \alpha \, (\eta, \chi).
\end{eqnarray}
With a choice of $u_{0 h}$ and $u_{1 h}$ as $L^2$-projection of $u_0$ and $u_1$, respectively, i.e., 
$$
(e'(0), \chi) = 0, \quad \mbox{and} \quad (e(0), \chi) = 0.
$$
Choosing $\chi = \theta + \epsilon \hat{\theta}$ in (\ref{En3.26}), we note that $\epsilon (\theta',\hat{\theta}) = \epsilon \frac{d}{dt} (\theta,\hat{\theta}) - \frac{\epsilon}{2} \|\theta\|^2.$ Setting corresponding discrete energy 
$${\cal E}_h^{(0)} (\theta)(t)  := \frac{1}{2} \big( \|\theta(t)\|^2 + \|A_h \hat{\theta}(t)\|^2\big)$$
with its extended energy
\begin{equation*}
{\cal E}^{(0)}_{h,\epsilon}(\theta)(t) = {\cal E}_h^{(0)}(\theta) (t)+ \epsilon (\theta, \hat{\theta}) 
+ \frac{1}{2} \alpha \, \epsilon \|\hat{\theta}(t)\|^2,
\end{equation*}
and 
\begin{eqnarray*}
F_0(t)&:=& (\alpha - \epsilon) \|\theta (t)\|^2 + {\epsilon} \|A^{1/2}_h \hat{\theta}\|^2\\
&=& 2 \epsilon {\cal E}_h^{(0)}(\theta) (t) + \left( \alpha - 2 \epsilon \right) \|\theta\|^2,
\end{eqnarray*}
to arrive  following similar to the proof of the Lemma~\ref{Lem:3.1}  at 
\begin{equation*} \label{En3.27}
\frac{d}{d t} {\cal E}^{(0)}_{h,\epsilon} (\theta)(t)  + \frac{4}{3} \epsilon  (1-\delta_0) \;{\cal E}_{h,\epsilon}^0 (\theta)(t) \leq \frac{1}{2\delta_0} \left( \frac{2}{\alpha} + \frac{\alpha}{2\lambda_1} \right) \;\left(\|\eta'(t)\|^2 +  \|\eta(t)\|^2\right).
\end{equation*}
Then, again, proceed in a similar to the lines of proof of the Lemma~\ref{Lem:3.1} to obtain
\begin{eqnarray*}
{\cal E}_h^{(0)} (\theta)(t) \leq 3 e^{-2 \delta t (1-\delta_0)} {\cal E}_h^{(0)}(\theta)(0) + \frac{1}{\delta_0} \left( \frac{2}{\alpha} + \frac{\alpha}{2\lambda_1} \right)  \int_0^t e^{-2 \delta  (1-\delta_0) (t - s)} \left(\|\eta'\|^2 + \|\eta\|^2 \right) \, ds,
\end{eqnarray*}
Since ${\cal E}_h^{(0)}(\theta)(0)=\frac{1}{2} \|\theta(0)\|$  and $\|\theta(0)\| \leq C\;h^2 \|u_0\|_{2}$, a use of estimates of 
$\|\eta\|$ and $\|\eta'\|$ with \eqref{estimate-e-t} and triangle inequality concludes the rest of the proof. \hfill $\blacksquare$


\section{Some Generalizations}
\setcounter{equation}{0}
In this section, we discuss some generalizations of our results to weakly damped wave equation with non-homogeneous forcing function, space dependent damping coefficient, viscous damping and compensation and weakly damped beam  equations.

\subsection{Inhomogeneous equations}
This subsection is on the weakly damped wave equation with non-homogeneous forcing function in abstract form as
\bea 
u'' + \alpha \, u' +A u = f, \;  t > 0 \nonumber 
\eea
with initial conditions
\bea
u(0) = u_0, \quad u'( 0) = u_1. \nonumber 
\eea
Here, $f = f(t) \in L^2.$ 
\begin{Teor}
Let $u_{\infty}$ be the unique solution of 
\begin{equation*} \label{EI.4}
A u_{\infty} = f, \quad \mbox{with} \quad u_{\infty} = 0 \quad \mbox{on} \quad \partial \Omega.
\end{equation*}
Then with $w(t) = u(t) - u_{\infty},$  there holds
\begin{equation*} \label{EI.5}
{\cal E}^{(j)}(w)(t) \leq 3 e^{- \delta t} {\cal E}^{(j)}(w)(0) = 3 e^{- \delta t} \left(\|w^{(j)}(0)\|^2 + \|A^{1/2} w^{(j-1)}(0)\|^2 \right).
\end{equation*}
Here, for $j = 1,$ there holds $w^{(1)} = w(0)$, and for $j > 1$, it follows that $w^{(j)}(0) = u^{(j)}(0).$
\end{Teor}
\noindent
{\bf Proof.} Now $w(t)$ satisfies 
\begin{eqnarray*} \label{EI.6}
w'' + \alpha \, w' + A w &=& 0, \;t\in (0, \infty), \\ \nonumber
w(0) = u_0 - u_{\infty}, &&  w'(0) = u_1.
\end{eqnarray*}
On following the technique for proving decay properties in Theorem 2.2 and Theorem 2.3, we complete the rest of the proof. \hfill $\blacksquare$

\begin{Rem} Following Theorem 2.3 and Theorem 2.4, we again arrive at
\begin{equation*}
{\cal E}^{(j)}_A(w)(t) \leq 3 e^{-2 \delta t} \, {\cal E}^{(j)}_A(w)(0). \label{new5.7}
\end{equation*}
\end{Rem}
\noindent
Thus, as in Remark 2.1, we find for $d=2$
\begin{equation*}
\|w(t)\|_{L^{\infty}} \leq C e^{-\delta t} \left(\|u_1\|_1 + \|u_0 - u_{\infty}\|_{2} \right). \label{new5.8}
\end{equation*}
This implies $u(t) \rightarrow u_{\infty}$ in $L^{\infty}(\Omega)$ as $t \rightarrow \infty$. \\\\
As in section 3, similar results holds  for the semidiscrete solution $w_h(t) = u_h(t) - u_{\infty, h}$ and hence, for some $1>\delta_0>0$, there holds $\|u_h(t)-u_{\infty, h}\|_{\infty} = O \left(e^{-\delta t(1-\delta_0)} \right)$. 
\begin{Rem}
In case $f(t) = O \left(e^{-\delta_0 t} \right)$, then also the solution decay exponentially with decay rate $\delta^* = \min \; (\delta_0, \delta)$.
\end{Rem}

\subsection{On space dependent damping term}
This subsection briefly focuses on the weakly damped wave equation with space dependent damping coefficient of the form, (see, \cite{chen-1}, \cite{CZ} and \cite{Rau}):
\bea 
u'' + \alpha \, u' + A u = 0,  \; t > 0  \nonumber 
\eea
with initial conditions
\bea
u(0) = u_0, \quad u'(0) = u_1. \nonumber 
\eea
Here, the space dependent damping coefficient $\alpha\in C^0( \bar{\Omega})$ satisfies 
$$
0 < \displaystyle{\min_{x \in \bar{\Omega}}} \, \alpha(x) = \alpha_1 \leq \alpha(x) \leq \alpha_2 = \max_{x \in \bar{\Omega}}\, \alpha(x).
$$

To indicate the decay property, for simplicity, assume that  $\alpha_1 \alpha_2 \leq \lambda_1$, where $\lambda_1$ is principal eigenvalue of the operator $A.$ An appropriate modification of the analysis of   Rauch \cite{Rau} shows that the continuous energy 
$$
{\cal E}^{(1)}(u)(t) = \frac{1}{2} \left(\|u'(t)\|^2 + \|A^{1/2} u(t)\|^2 \right),
$$
decays like 
\begin{equation*}
{\cal E}^{(1)}(u)(t) \leq \max \left(4, \frac{\alpha_1^2}{2 \lambda_1} \right) e^{-\alpha_1 t} {\cal E}^{(1)}(u)(0). \label{New5.12}
\end{equation*}
Similarly, by differentiating  $j$ times in the temporal variable, it follows easily that 
\begin{equation*}
{\cal E}^{(j)}(u)(t) \leq \max \left(4, \frac{\alpha_1^2}{2 \lambda_1} \right) e^{-\alpha_1 t} {\cal E}^{(j)}(u)(0). \label{New5.13}
\end{equation*}

\noindent
For the corresponding semidiscrete system: Find $u_h(t) \in S_h^0$ such that 
\begin{equation}
(u''_h, \chi_h) + (\alpha u'_h, \chi_h) + a (u_h, \chi_h) = 0 \;\;\; \forall \chi_h \in S_h^0. \label{New5.14}
\end{equation}

\noindent
Setting  $w_h = e^{(\alpha_1 /2) t} u_h(t),$ we now rewrite \eqref{New5.14} in terms of $w_h$ as
\begin{equation}
(w''_h, \chi_h)  + a (w_h, \chi_h) +(  (\frac{\alpha_1^2}{4}-\frac{\alpha \alpha_1}{2}) w_h,\chi_h)+( (\alpha-\alpha_1) w'_h, \chi_h) = 0 \;\;\; \forall \chi_h \in S_h^0. \label{New5.14-1}
\end{equation}
Now choose $\chi_h=w_h'$ in \eqref{New5.14-1} and define 
$$
{\cal I}_h (w_h)(t) = {\cal E}^{(1)}_h (w_h)(t) + \frac{1}{2} \int_{\Omega} \left(\frac{\alpha_1^2}{4} - \frac{\alpha \alpha_1}{2} \right)\, |w_h|^2 \, dx.
$$
Then, as $(\alpha-\alpha_1)\geq 0$, there holds 
\begin{equation*}\label{modified-energy}
\frac{d}{dt} {\cal I}_h (w_h)(t) = -( (\alpha-\alpha_1) w'_h, w_h') \leq 0,
\end{equation*}
and an integration with respect to time shows
\begin{equation}\label{modified-energy-1}
{\cal I}_h (w_h)(t) \leq {\cal I}_h (w_h)(0).
\end{equation}
Note that $\frac{\alpha_1^2}{4} - \frac{\alpha \alpha_1}{2} \leq - \frac{\alpha_1^2}{4} < 0.$
Since $e^{\frac{\alpha_1 }{2}t} u_h'=w_h'-\frac{\alpha_1}{2} w_h(t),$ it follows using $(a-b)^2 \leq 2(a^2 + b^2)$ that 
\begin{eqnarray*}
e^{\alpha_1 t} {\cal E}_h^{(1)}(u_h)(t) &=& \frac{1}{2} \left(\|\big(w'_h - \frac{\alpha_1}{2} w_h \big) (t) \|^2 + \|A^{1/2}_h w_h(t)\|^2 \right) \\
&\leq& \left(\|w'_h(t)\|^2 + \frac{\alpha_1^2}{4} \|w_h(t)\|^2 + \frac{1}{2} \|A^{1/2}_h w_h(t)\|^2 \right) \\
&\leq& 2 {\cal E}_h^{(1)}(w_h)(t) + \frac{\alpha_1^2}{4} \|w_h(t)\|^2 - \frac{1}{2} \|A^{1/2}_h w_h(t)\|^2.
\end{eqnarray*}
Since $ \alpha_1 \alpha_2 + 2 \left(\frac{\alpha_1^2}{4} - \frac{\alpha \alpha_1}{2} \right) \geq \frac{\alpha_1^2}{2}$, we obtain
$$
e^{\alpha_1 t} {\cal E}_h^{(1)}(u_h)(t) \leq 2 \, {\cal I}_h (w_h)(t) + \frac{1}{2} \alpha_2 \alpha_1 \|w_h(t)\|^2 - \frac{1}{2} \|A^{1/2}_h w_h(t)\|^2.
$$
A use of the Poincar{\'e} inequality $\|w_h(t)\|^2 \leq \frac{1}{\lambda_1} \|A^{1/2}_h w_h(t)\|^2$ shows
$$
e^{\alpha_1 t} {\cal E}_h^{(1)}(u_h)(t) \leq 2 \, {\cal I}_h(w_h(t)) +  \frac{1}{2} \left(\alpha_1 \alpha_2 - \lambda_1 \right) \| w_h(t)\|^2.
$$
If $\alpha_1^2 < \alpha_1 \alpha_2 \leq \frac{\lambda_1}{2}$, then $\alpha_1 \alpha_2 - \lambda_1\leq 0$. Thus, a use of \eqref{modified-energy-1} yields
\begin{equation}
 {\cal E}_h^{(1)}(u_h)(t) \leq 2 \,e^{-\alpha_1 t}\, {\cal I}_h (w_h)(t) \leq 2 \,e^{-\alpha_1 t}\, {\cal I}_h (w_h)(0) . \label{New5.15}
\end{equation}
Since $ \left(\frac{\alpha_1^2}{4} - \frac{\alpha \alpha_1}{2} \right)  \leq \frac{\alpha_1^2}{4},$ we note with $\dsp{e^{-\frac{\alpha_1}{2} t} w_h'(t)=
u_h'(t) + \frac{\alpha_1}{2} u_h(t)}$ and the Poincar{\'e} inequality that
\begin{eqnarray}
2 {\cal I}_h (w_h(0)) &=&  \|u_h'(0) + \alpha u_h(0)\|^2 
+ \|A_h^{1/2} u_h(0)\|^2+ \int_{\Omega} \left(\frac{\alpha_1^2}{4} - \frac{\alpha \alpha_1}{2} \right) |u_h(0)|^2 \;dx\nonumber \\
&\leq& 2\;{\cal E}^{(1)}_h (u_h)(0) +  \frac{1}{4} \alpha_1^2  \|u_h(0)\|^2\;dx \leq 2\;{\cal E}^{(1)}_h (u_h)(0)
+\frac{1}{4\lambda_1} \alpha_1^2  \|A_h^{1/2} u_h(0)\|^2  \nonumber \\
&\leq& \max\left(2, \frac{\alpha_1^2}{4 \lambda_1} \right) {\cal E}^{(1)}_h (u_h)(0). \label{New5.17}
\end{eqnarray}
On substitution of 
(\ref{New5.17}) in (\ref{New5.15}), we arrive at
$$
{\cal E}^{(1)}_h (u_h)(t) \leq \max\left(2, \frac{\alpha_1^2}{4 \lambda_1} \right) e^{-\alpha_1 t} {\cal E}^{(1)}_h (u_h)(0).
$$
Similarly,
$$
{\cal E}^{(j)}_h (u_h)(t) \leq \max  \left(2, \frac{\alpha_1^2}{4 \lambda_1} \right) e^{-\alpha_1 t} {\cal E}^{(j)}_h (u_h)(0).
$$
Moreover, we derive all the error estimates as in Section 3. In particular, when $d=2$ and for small $0<\delta_0<1$, there holds
$$
\|(u - u_h)(t)\|_{\infty} \leq C \left( \log \big(\frac{1}{h}\big) \right) h^2 \sqrt{t} e^{- \frac{1}{2} \alpha_1(1-\delta_0) t}.
$$
\begin{Rem}
In section 3, since $\alpha$ is a constant, the decay rate  is $O\left(e^{- \frac{\alpha}{2}(1-\delta_0) t} \right)$ for small $\delta_0>0$, provided
 $\alpha_1\alpha_2=\alpha^2 \leq \lambda_1.$ In fact, the analysis of this subsection  improves the decay  rate compared to the decay rate  in the Section 3.
\end{Rem}

\subsection{On viscous damping and compensation}
This subsection is  on improved decay rates due to both viscous damping and compensation, which is influenced by Chen \cite{chen-3}.

Now, consider  the wave equation with  positive constant  viscous damping and compensation terms which is written in abstract form as: 
\bea 
u'' + \alpha \, u' + \beta\; u + A u  = 0, \;  t > 0, \nonumber 
\eea
with initial conditions
\bea
u(0) = u_0, \quad u_t(0) = u_1. \nonumber 
\eea
Here, $\alpha$ and $\beta$ are called the viscous damping and compensation coefficient, respectively. When $A = -\Delta$, this problem was discussed in \cite{chen-3}, and improved exponential decay rates were established. For a general second order linear  self-adjoint positive elliptic operator, appropriate modification provides the following improved decay estimates for the energy. 
\begin{Teor}
For any $\delta >0$ with
\begin{equation}\label{condition:delta}
\alpha = \delta (3+\delta)\;\;\mbox{and}\;\; \beta= \delta ( 2 + 3 \delta + 2 \delta^2),
\end{equation}
the energy 
$$
{\cal E}^{(1)}(u)(t) = \frac{1}{2} \left(\|u'(t)\|^2 + \|A^{1/2} u(t)\|^2 \right),
$$
decays exponentially, that is,
\begin{equation}
{\cal E}^{(1)}(u)(t) \leq C(\lambda_1,\delta) \; e^{- \delta t} {\cal E}^{(1)}(u)(0), \label{vdc:1}
\end{equation}
where the positive constant $C=O(\delta^3).$
\end{Teor}
Note that for large $\alpha$ and $\beta$, it is possible to derive decay rate $\delta>0$, which remains large. 
Moreover, for a given $\delta >0$ with \eqref{condition:delta},
$u_0\in D(A^{(k)})$ and $u_1\in D(A^{(k-1/4)})$ for $k>1,$  there holds using the arguments to arrive at \eqref{vdc:1} and using induction
$${\cal E}_{A^{(k)}}^{(j)}(u)(t) \leq C(\lambda_1; \delta)\,e^{-\delta t} \, {\cal E}_{A^{(k)}}^{(j)}(u)(0),$$
where ${\cal E}_{A^{(k)}}^{(j)}(u)(t)=\frac{1}{2} \left(\|A^{(k-1/4)} u^{(j)}(t)\|^2 + \|A^{(k)} u^{(j-1)}(t)\|^2 \right).$

Now, the corresponding semidiscrete system is to seek $u_h(t) \in S_h^0$ such that 
\begin{equation}
(u''_h, \chi_h) + \alpha (u'_h, \chi_h) + a (u_h, \chi_h) + \beta (u_h, \chi_h) = 0 \;\;\; \forall \chi_h \in S_h^0. \label{vdc:2}
\end{equation}
With a choice of $\chi_h= u_h' + \delta u_h$ in \eqref{vdc:2}, it follows using definition $A_h$ as in \eqref{discrete-A} with the energy
$$
{\cal E}^{(1)}_h (u_h)(t) = \frac{1}{2} \left(\|u_h'(t)\|^2 + \|A^{1/2}_h u_h(t)\|^2 \right),
$$
and extended energy
$$
{\cal E}^{(1)}_{\delta,h}(u_h)(t) =  {\cal E}^{(1)}_h (u_h)(t)  +\frac{1}{2} (\beta + \delta \alpha) \|u_h(t)\|^2 + \delta (u_h', u_h),
$$
that 
\begin{equation}\label{vdc:3}
\frac{d}{dt} {\cal E}^{(1)}_{\delta,h}(u_h)(t) + F_h(t)= 0,
\end{equation}
where
$$F_h(t) :=  (\alpha-\delta) \|u_h'(t)\|^2 + \delta \|A_h^{1/2} u_h(t)\|^2 + \beta \delta \|u_h(t)\|^2.$$
Since from \eqref{condition:delta}, the condition
\begin{equation*}\label{condition:db}
\beta \delta \geq \frac{\delta}{2} ( \delta + \beta + \alpha \delta ),
\end{equation*}
shows using  $ -\delta^2 ( u_h', u_h) \geq -\big(  (\delta^2/2) \|u_h'\|^2 + (\delta^2/2) \|u_h\|^2\big)$  that
\begin{eqnarray}\label{vdc:4}
F_h(t) \geq \delta {\cal E}^{(1)}_{\delta,h}(u_h)(t)  + \frac{\delta}{2}  (3+ \delta) \|u_h\|^2 \geq \delta {\cal E}^{(1)}_{\delta,h}(u_h)(t).
\end{eqnarray}
On substitution of \eqref{vdc:4} in  \eqref{vdc:3}, we arrive at 
$$ 
\frac{d}{dt} {\cal E}^{(1)}_{\delta,h}(u_h)(t) + \delta {\cal E}^{(1)}_{\delta,h}(u_h)(t) \leq \frac{d}{dt} {\cal E}^{(1)}_{\delta,h}(u_h)(t) + F_h(t)= 0,
$$
and hence, an integration with respect to time yields
\begin{equation}\label{vdc:5}
 {\cal E}^{(1)}_{\delta,h}(u_h)(t) \leq e^{-\delta t} \; {\cal E}^{(1)}_{\delta,h}(u_h)(0).
 \end{equation}
 Again a use of \eqref{condition:delta} shows
 \begin{eqnarray}\label{vdc:6}
 {\cal E}^{(1)}_{\delta,h}(u_h)(t) &\geq& \frac{1}{2} \Big( \|u_h'(t)\|^2 + \|A_h^{1/2} u_h(t)\|^2 + (\beta + \delta \alpha) \|u_h(t)\|^2\Big) 
 -\frac{1}{4} \|u_h'(t)\|^2 - \delta^2 \|u_h(t)\|^2\nonumber\\
 &=& \frac{1}{4} \Big( \|u_h'(t)\|^2 +2 \|A_h^{1/2} u_h(t)\|^2 \Big) + \big( \frac{1}{2} (\beta + \delta \alpha)-\delta^2\big)  \|u_h(t)\|^2\nonumber\\
 &\geq & \frac{1}{4} \Big( \|u_h'(t)\|^2 + \|A_h^{1/2} u_h(t)\|^2 \Big)  + \frac{1}{2} \big(2\delta + 4 \delta^2 + 3 \delta^3\big) \|u_h(t)\|^2\nonumber\\
 &\geq& \frac{1}{2} {\cal E}^{(1)}_h(u_h) (t). 
 \end{eqnarray}
 For obtaining an upper bound, we note using  \eqref{condition:delta}, $\delta (u'_h,u_h) \leq (1/2) ( \|u_h'\|^2 + \delta^2 \|u_h\|^2)$ and Poincar{\'e} inequality \eqref{Ph1} 
 \begin{eqnarray}\label{vdc:7}
 {\cal E}^{(1)}_{\delta,h}(u_h)(t) &\leq& \frac{1}{2} \Big( 2 \|u_h'(t)\|^2 + \|A_h^{1/2} u_h(t)\|^2 + (\beta + \delta \alpha + \delta^2) \|u_h(t)\|^2 \Big)  
  \nonumber\\
 &\leq&  \|u_h'(t)\|^2 + \frac{1}{2 \lambda_1} (1 +\beta + \delta \alpha + \delta^2) \|A_h^{1/2} u_h(t)\|^2 \Big)  
  \nonumber\\ 
  &\leq& \frac{1}{2 \lambda_1} \big( 2 \lambda_1 + (2\delta + 7 \delta^2 + 3 \delta^3\big)\; {\cal E}^{(1)}_h(u_h)(t).
  \end{eqnarray}
 With $ \frac{1}{2} C(\lambda_1, \delta) = \frac{1}{2 \lambda_1} \Big( 2 \lambda_1 + (2\delta + 7 \delta^2 + 3 \delta^3)\Big)= O(\delta^3)$, we arrive  from  \eqref{vdc:6}-\eqref{vdc:7}   at
 \begin{equation}\label{equiv:energy}
 \frac{1}{2}  {\cal E}^{(1)}_h(u_h)(t) \leq  {\cal E}^{(1)}_{\delta,h}(u_h)(t) \leq \frac{1}{2} \;C(\lambda_1,\delta) \;{\cal E}^{(1)}_h(u_h)(t) .
 \end{equation}
 On substitution in \eqref{vdc:5}, we obtain
 \begin{equation}\label{vdc:8}
 {\cal E}^{(1)}_h(u_h) (t)  \leq C(\lambda_1,\delta)\;e^{-\delta t} \;{\cal E}^{(1)}_h(u_h) (0). 
 \end{equation}
 Moreover, following the similar line of arguments, there holds for $j\geq 1$
  \begin{equation*}
 {\cal E}^{(j)}_h(u_h) (t)  \leq C(\lambda_1,\delta)\;e^{-\delta t} \;{\cal E}^{(j)}_h(u_h)(0). 
 \end{equation*}
Further, a use of definition of $A_h$ in \eqref{discrete-A} yields  
 \begin{equation*}
 {\cal E}^{(j)}_{A_h}(u_h) (t)  \leq C(\lambda_1,\delta)\;e^{-\delta t} \;{\cal E}^{(j)}_{A_h}(u_h) (0). 
 \end{equation*}
 
 Following the argument that leads to \eqref{vdc:8} and also the error analysis in section 3, the following  optimal error estimates for $\delta >0$ and for any small $\delta_0>0$ hold:
 $$
\|u(t) - u_h(t)\| + h \|\nabla(u(t)-u_h(t)\| \leq C h^2 \,(1+t)^{1/2} \, e^{-\frac{\delta}{2} (1-\delta_0) t},
$$
 and for $d=2$
 $$\|u(t) - u_h(t)\|_{L^{\infty}} \leq C h^2 \left( \log \left(\frac{1}{h} \right) \right)(1+ t)^{1/2} e^{- \frac{\delta}{2} (1-\delta_0) t}.$$
 
\subsection{ On  weakly damped beam  equations} 
This subsection is on the beam equation with a weakly damping term, see \cite{GR}.

For a convex polygonal or polyhedral domain  $\Omega$  in $\mathbb{R}^d$ with boundary $\partial \Omega$  and   fixed positive constant
$\alpha$, the problem is to find $u(x,t)$ for $(x,t) \in \Omega \times (0, \infty)$ satisfying 
\bea 
u_{tt} + \alpha \, u_t + \Delta^2 u = 0, \; x \in \Omega, \; t > 0, \label{B1.1}
\eea
with initial conditions
\bea
u(x, 0) = u_0(x), \quad u_t(x, 0) = u_1(x), \; x \in \Omega, \label{B1.2}
\eea
and homogeneous clamped  boundary conditions
\bea
u = \frac{\partial u}{\partial \nu}= 0, \quad (x, t) \in \partial \Omega \times (0, \infty), \label{B1.3}
\eea
where $\nu$ is the outward unit normal to the boundary $\partial \Omega.$ 

With $A= \Delta^2$  and $D(A)= H^4(\Omega) \cap H^2_0(\Omega)$, results of the  previous sections  remain valid in the present case with appropriate changes. For semidiscrete FEM, choose $S_h^0$ be a finite element subspace of $H_0^2(\Omega)$ satisfying the following approximation property:
$$
\inf_{\chi \in S_h^0} \sum_{j = 0}^2 h^j \|v - \chi\|_{H^j(\Omega)} \leq C h^3  |v|_{H^3(\Omega)}.
$$ 

Then, the rest of the decay property holds similarly. Based on the arguments in Section 3, for small $\delta_0>0$ the following estimates are easy to hold
$$
\|(u - u_h)(t)\|_j = O \left(h^{3 - j} e^{-\delta (1-\delta_0) \, t} \right), \; j = 1, 2.
$$

Instead of homogeneous clamped boundary conditions, we can use either hinged boundary conditions or simply supported boundary conditions, and with appropriate modifications, similar results can be derived. 

\section{Numerical Experiments}
\setcounter{equation}{0}
This section focusses on some numerical experiments, whose results confirm our theoretical findings in Sections 3 and 4.

\subsection{Completely Discrete Scheme}
Let $k > 0$ be the time step and let $t_n = n k, \; n \geq 0$. Set 
$\varphi^n = \varphi(t_n)$,
$$
{\bar \partial}_t \varphi^n = \frac{\varphi^n - \varphi^{n-1}}{k} \quad \mbox{and} \quad \partial_t \varphi^n = \frac{\varphi^{n + 1} - \varphi^{n}}{k}
$$
with ${\bar \partial}^0_t \varphi^n = \varphi^n$.  Define  
$$
{\bar \partial}^{(j + 1)}_t \varphi^n = \frac{1}{k} \left( {\bar \partial}^{j}_t \varphi^n - 
{\bar \partial}^{j}_t \varphi^{n - 1} \right), \;\; j \geq 0
$$
and  $\varphi^{n + \frac{1}{2}} = \frac{\varphi^{n + 1} + \varphi^n}{2}.$  Set
\begin{eqnarray*}
\delta_t \varphi^n &=&  \frac{\varphi^{n+1} - \varphi^{n-1}}{2 k} = {\bar \partial}_t \varphi^{n + \frac{1}{2}} = \frac{\varphi^{n + \frac{1}{2}} - \varphi^{n - \frac{1}{2}}}{k}, \\
{\hat \varphi}^n &=& \frac{1}{4} \left(\varphi^{n+1} + 2 \varphi^n + \varphi^{n-1} \right) = \frac{1}{2} \left(\varphi^{n + \frac{1}{2}} + \varphi^{n - \frac{1}{2}} \right), \\
\partial_t {\bar \partial}_t \varphi^n &=& \frac{1}{k^2} \left(\varphi^{n + 1} - 2 \varphi^n + \varphi^{n -1} \right) = \frac{1}{2 k} \left(\varphi^{n + \frac{1}{2}} - \varphi^{n - \frac{1}{2}} \right) = \frac{1}{k} \left(\partial_t \varphi^n - {\bar \partial}_t \varphi^n \right).
\end{eqnarray*}
The discrete time finite element approximations $U^n$ of
$u(t_n)$ is defined as solution of 
\begin{eqnarray}
(\partial_t {\bar \partial}_t U^n, \chi) + \alpha (\delta_t U^n, \chi) 
+a ( {\hat U}^n, \chi) = 0, \; \chi \in S_h^0, \; n \geq 1  \label{Feqn3.1}
\end{eqnarray}
with $U^0 = u_{0,h}$ and $U^1 = u_{1,h},$
where $u_{0,h}, u_{1,h} \in S_h^0$ are appropriate approximations to be
defined later. 

We now define the discrete energy
\begin{eqnarray*}
{\cal E}^n(U) = \frac{1}{2} \left( \| \partial_t U^n \|^2 + \|A^{1/2} U^{n + \frac{1}{2}} \|^2 \right), \; n \geq 0. \label{Feqn3.2}
\end{eqnarray*}
Choose   $\chi = \delta_t U^n$ in (\ref{Feqn3.1}) to obtain
\begin{equation}
\left(\partial_t {\bar \partial}_t U^n, \delta_t U^n \right) + \alpha \|\delta_t U^n\|^2 + a (  {\hat U}^n,  \delta_t U^n ) = 0. \label{Feqn3.3}
\end{equation}
Note that 
\begin{equation}
\left(\partial_t {\bar \partial}_t U^n, \delta_t U^n \right) = \frac{1}{2k} \left(\partial_t U^n - {\bar \partial}_t U^n, \partial_t U^n + {\bar \partial}_t U^n \right) = \frac{1}{2k} \left(\|\partial_t U^n\|^2 - \|\partial_t U^{n-1}\|^2 \right), \label{Feqn3.4}
\end{equation}
and
\begin{eqnarray}
a( {\hat U}^n,\delta_t U^n ) &=& \frac{1}{2k} a( (U^{n + \frac{1}{2}} + U^{n - \frac{1}{2}} ),  (U^{n + \frac{1}{2}} - U^{n - \frac{1}{2}} ) ) \nonumber \\
&=& \frac{1}{2k} \left(\|A^{1/2} U^{n + \frac{1}{2}} \|^2 - \|A^{1/2} U^{n - \frac{1}{2}} \|^2 \right). \label{Feqn3.5} 
\end{eqnarray}
Substituting (\ref{Feqn3.4})-(\ref{Feqn3.5}) in (\ref{Feqn3.3}), we obtain
\begin{equation*}
{\cal E}^n(U) - {\cal E}^{n-1}(U) + \alpha \, k \|\delta_t U^n\|^2 = 0. 
\end{equation*}
Taking summation for $n = 1$ to $m$, we arrive at
\begin{equation*}
{\cal E}^m(U) + \alpha \, k \sum_{n = 1}^m \|\delta_t U^n\|^2 = {\cal E}^0(U). \label{Feqn3.6} 
\end{equation*}
Therefore, the discrete energy satisfies
$$
{\cal E}^n(U) \leq {\cal E}^0(U).
$$
For numerical experiments, examples 1, 2, and 6 are related to the homogeneous weakly damped wave equation with various damping parameter values.  Examples 3, 4, 5, and 7 are related to the weakly damped wave equation with a nonhomogeneous forcing function. In examples 1-7, the equations are solved up to the final time $T = 1.0$ with the time step $k = h^2$. The numerical experiments are performed using FreeFem++ with piecewise linear elements \cite{ff++}.

In each case, the experimental convergence rate of the error is computed using
\[
\mbox{Rate} = \frac{\log(E_{h_{i}})-\log (E_{h_{i+1}})}{\log (\frac{h_i}{h_{i+1}})},
\]
where $E_{h_{i}}$ denotes the norm of the error using $h_i$ as the spatial discretizing parameter at $i^{\small{\mbox{th}}}$ stage. \\

\noindent
{\bf Example 1.} For the weakly damped wave equation:
\beas
u_{tt} + \alpha \, u_t - \Delta u = 0, \; (x_1, x_2) \in \Omega = (0, 1) \times (0, 1), \; t > 0
\eeas
with initial conditions
\beas
u(x_1, x_2, 0) = \sin (\pi x_1) \sin(\pi x_2), \quad u_t(x_1, x_2, 0) = \left(- \frac{\alpha}{2} + \frac{1}{2}\sqrt{\alpha^2 - 8 \, \pi^2} \right) \sin (\pi x_1) \sin(\pi x_2), \; (x_1, x_2) \in \Omega 
\eeas
and homogeneous Dirichlet boundary condition,
the exact solution is given by
$$
u(x_1, x_2, t) = e^{\left(- \frac{\alpha}{2} + \frac{1}{2}\sqrt{\alpha^2 - 8 \, \pi^2} \right) t} \sin (\pi x_1) \sin(\pi x_2).
$$
Table \ref{ex-t1} shows the errors and rate of converges in $L^2, \; L^{\infty}$ and $H^1$-norms, confirming our theoretical findings. \\
\begin{table}[!h]
$$
\begin{array}{|c|cccccc|} \hline
&   & & \alpha = 0.9  &  &  &  \\ \cline{2-7}
N &  \|u-u_{h}\| & \mbox{Rate} & \|u-u_{h}\|_{\infty} &\mbox{Rate}
& \|u-u_{h}\|_1 &
\mbox{Rate} \\ \hline
6  & 1.22937(-3) &   -    & 9.09229(-4)& -      &    1.20156(-2) & -        \\ 
12 & 4.5851(-4) & 1.3955 & 3.61551(-4)& 1.30483 & 4.45790(-3) & 1.40294\\
18 & 2.06364(-4) & 2.23398 & 1.63361(-4)& 2.22307  & 3.10802(-3) & 1.00932  \\
24 & 1.17616(-4) & 1.96074 & 9.47278(-5) & 1.90055 & 2.33173(-3) & 1.00222  \\    
30 & 7.47123(-5) & 2.31882 &  5.95625(-5)& 2.37094 & 1.74253(-3) & 1.48841  \\ \hline
\end{array}
$$
\caption{Example 1: Errors and rate of convergences in $\|u-u_h\|, \; \|u - u_h\|_{\infty}$ and $\|u-u_h\|_1$.}
\label{ex-t1}
\end{table}

\begin{figure}[htb!]
\centering
\begin{subfigure}{.45\textwidth}
  \centering
  \includegraphics[width=.9\linewidth]{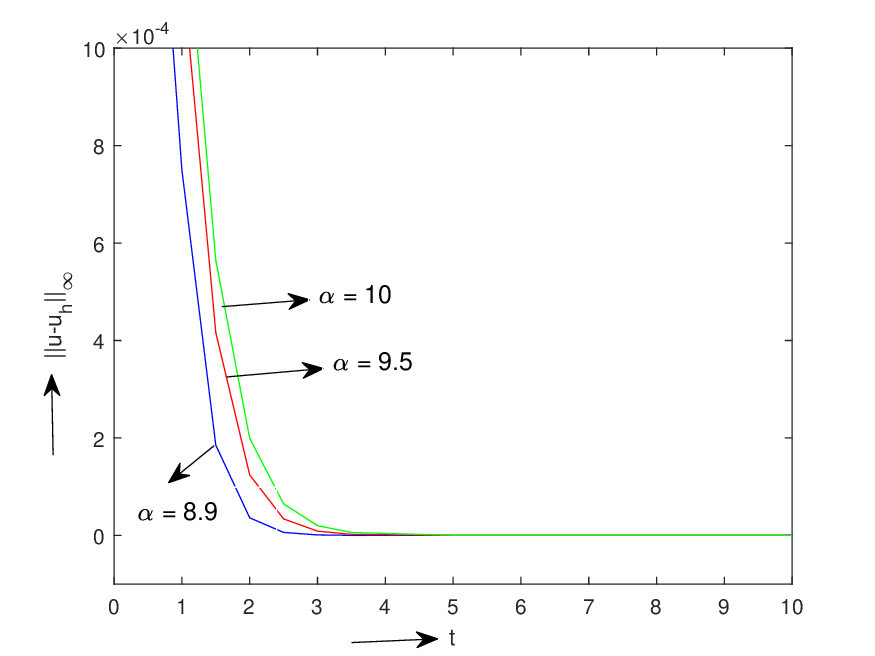}
  \caption{The decay estimate in $L^{\infty}$-norm.}
  \label{fig:Ex2.1}
\end{subfigure}%
\begin{subfigure}{.45\textwidth}
  \centering
  \includegraphics[width=.9\linewidth]{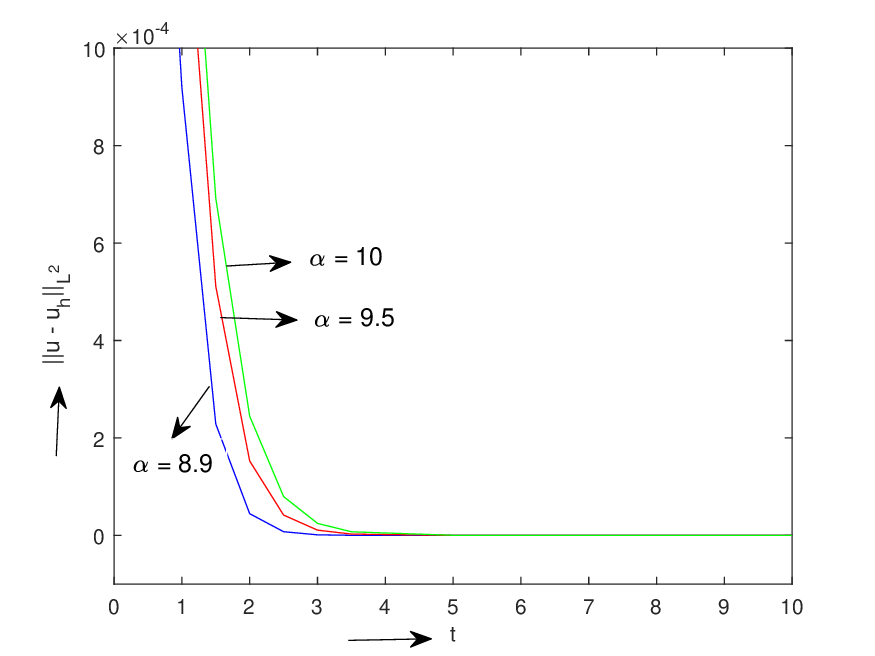}
  \caption{The decay estimate in $L^{2}$-norm.}
  \label{fig:Ex2.2}
\end{subfigure}
\begin{subfigure}{.45\textwidth}
  \centering
  \includegraphics[width=.9\linewidth]{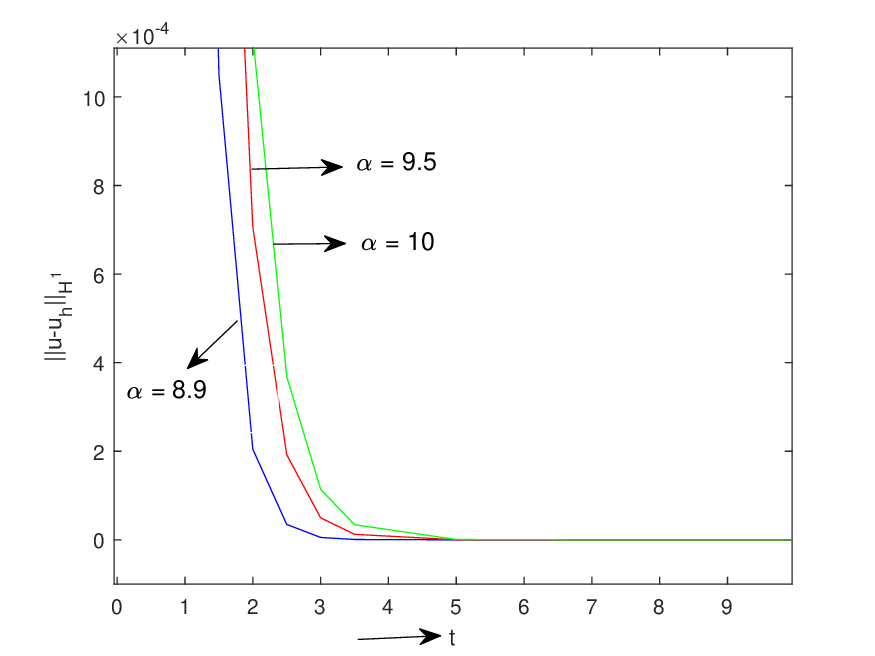}
  \caption{The decay estimate in $H^{1}$-norm.}
  \label{fig:Ex2.3}
\end{subfigure}
\caption{Example 1: The decay estimates in $L^{\infty}, \, L^2$ and $H^{1}$-norms.}
\label{ex-f1}
\end{figure}
\noindent
From Figure \ref{ex-f1}, we observe that $\alpha = 8.9$ exponentially decay faster than $\alpha = 9.5$ and $\alpha = 10$. This confirms that exponentially decay phenomenon for all the three norms $L^{\infty}, \, L^2$ and $H^{1}$. \\

\begin{figure}[!h]
  \centering
  \includegraphics[width=.5\linewidth]{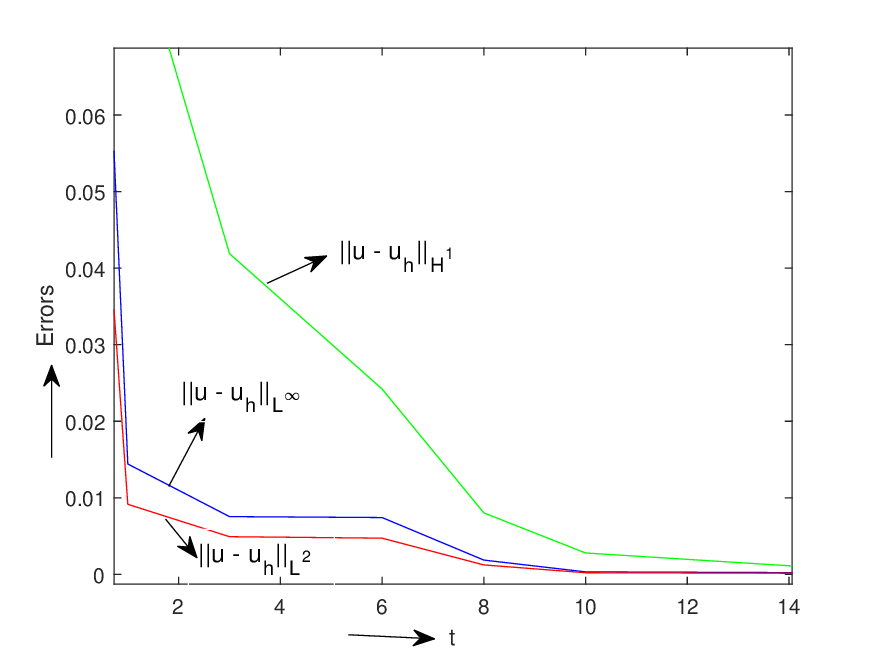}
  \caption{Example 2: The decay estimates in $L^{\infty}, \, L^2$ and $H^{1}$-norms for $\alpha_0 = 1$ and $\gamma = \frac{1}{2}$.}
  \label{ex-f2}
\end{figure}

\noindent
{\bf Example 2.} \cite{Rau} For the weakly damped wave equation with space dependent damping coefficient of the form
\beas
u_{tt} + \alpha(x_1,x_2) \, u_{t} - \Delta u = f(x_1, x_2, t), \; (x_1, x_2) \in \Omega = (1,2) \times (1, 2), \; t > 0
\eeas
with initial conditions
\beas
u(x_1, x_2, 0) = u_0(x_1, x_2), \quad u_t(x_1, x_2, 0) = u_1(x_1, x_2), \; (x_1, x_2) \in \Omega
\eeas
and homogeneous Dirichlet boundary condition, 
where $\alpha(x_1,x_2) = \alpha_0 (x_1^2+x_2^2)^{- \gamma/2}$ with some $\alpha_0 > 0$ and $\gamma = [0, 1)$, 
we compute the unknowns $f, \; u_0$ and $u_1$ with the help of the exact solution
$$
u(x_1, x_2, t) = e^{-\pi t} \sin (\pi x_1) \sin (\pi x_2).
$$
In Table \ref{ex-t2}, the errors and rate of converges in $L^2, \; L^{\infty}$ and $H^1$-norms are shown, and in Figure \ref{ex-f2}, we observe that errors decay exponentially.\\

\begin{table}[!h]
$$
\begin{array}{|c|cccccc|} \hline
&   & & \alpha_0 = 1, \gamma = \frac{1}{2}  &  &  &  \\ \cline{2-7}
N &  \|u-u_{h}\| & \mbox{Rate} & \|u-u_{h}\|_{\infty} &\mbox{Rate}
& \|u-u_{h}\|_1 &
\mbox{Rate} \\ \hline
6  & 3.87107(-2) &   -    &  5.97962(-2)& -      &    2.28391(-1) & -        \\ 
12 & 5.71754(-3) & 2.71788 &  8.9329(-3)&  2.70172  & 6.47137(-2) & 1.79208 \\
18 & 2.73750(-3) & 2.29648 & 4.27987(-3)&  2.29437  & 4.16378(-2) & 1.37498  \\
24 & 1.36909(-3) & 2.10169 &  2.16636(-3) & 2.06522  & 3.11371(-2) & 0.881475  \\    
30 & 8.88173(-4) & 1.93243 & 1.41087(-3)&   1.93243  & 2.47104(-2) & 1.03234  \\ \hline
\end{array}
$$
\caption{Example 2: Errors and rate of convergences in $\|u-u_h\|, \; \|u - u_h\|_{\infty}$ and $\|u-u_h\|_1$.}
\label{ex-t2}
\end{table}

\begin{figure}[htb!]
  \centering
  \includegraphics[width=.55\linewidth]{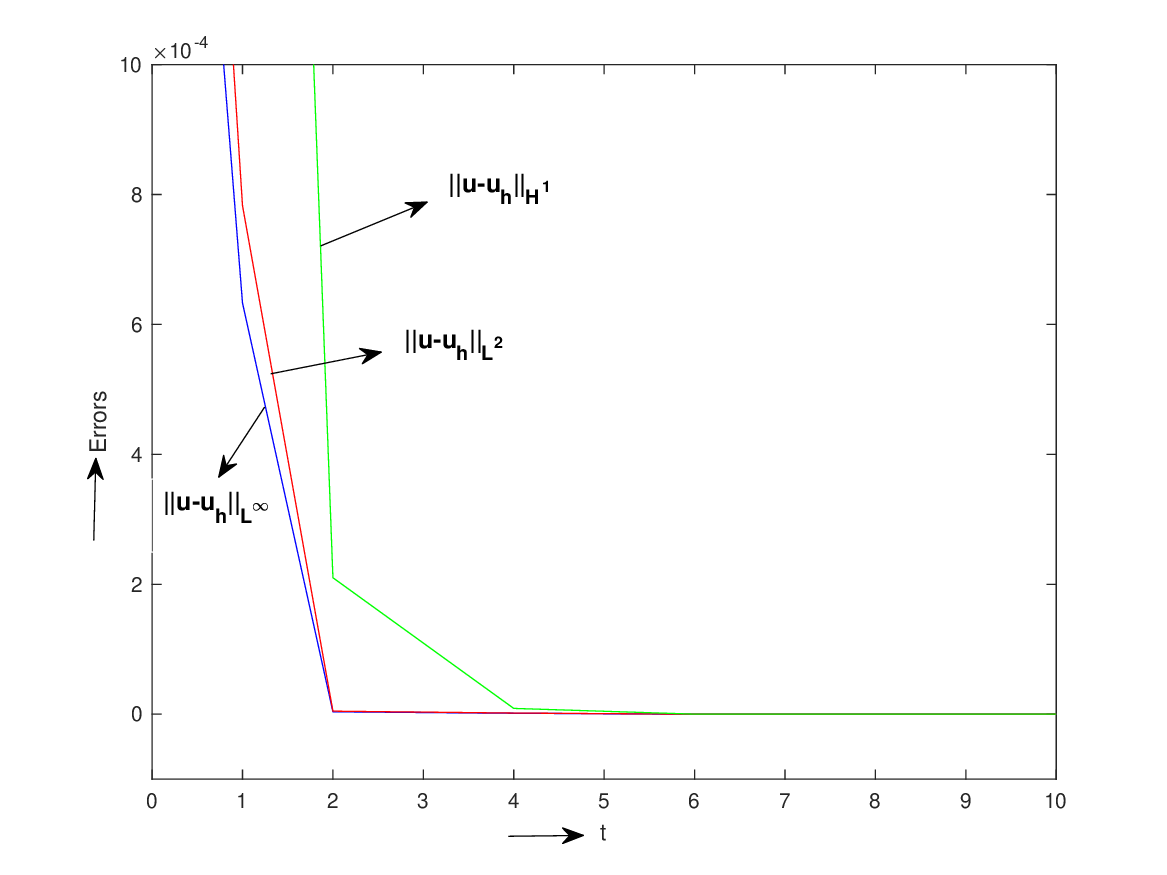}
  \caption{Example 3: The decay estimates in $L^{\infty}, \, L^2$ and $H^{1}$-norms for $\alpha = 4$.}
  \label{ex-f3}
\end{figure}

\noindent
{\bf Example 3.} For the semilinear weakly damped wave equation, see \cite{JMB} and \cite{Kap}
\beas
u_{tt} + \alpha \, u_{t} - \Delta u + f(u) = g(x_1, x_2, t), \; (x_1, x_2) \in \Omega = (0,1) \times (0, 1), \; t > 0
\eeas
with initial conditions
\beas
u(x_1, x_2, 0) = u_0(x_1, x_2), \quad u_t(x_1, x_2, 0) = u_1(x_1, x_2), \; (x_1, x_2) \in \Omega
\eeas
and homogeneous Dirichlet  boundary condition,
where $f(u) = u^3 - u,$ we compute the unknowns $g, \; u_0$ and $u_1$ with the help of the exact solution
$$
u(x_1, x_2, t) = e^{-\pi t} \sin (\pi x_1) \sin (\pi x_2).
$$
The errors and rate of converges in $L^2, \; L^{\infty}$ and $H^1$-norms are shown in the Table \ref{ex-t3}. In Figure \ref{ex-f3}, we observe that errors decay exponentially.

\begin{table}[!h]
$$
\begin{array}{|c|cccccc|} \hline
&   & & \alpha = 4  &  &  &  \\ \cline{2-7}
N &  \|u-u_{h}\| & \mbox{Rate} & \|u-u_{h}\|_{\infty} &\mbox{Rate}
& \|u-u_{h}\|_1 &
\mbox{Rate} \\ \hline
8  & 1.83112(-3) &   -    & 1.48149(-3)& -      &    2.56137(-2) & -        \\ 
16 & 6.92957(-4) & 1.34342 & 5.63024(-4)& 1.33757 & 1.04974(-2) & 1.23321 \\
24 & 3.25962(-4) & 1.84440 & 2.61394(-4)& 1.87645  & 6.27909(-3) & 1.25676  \\
32 & 1.8708(-4) & 1.89517 & 1.50115(-4) & 1.89306 & 4.45367(-3) & 1.17243  \\    
40 & 1.21071(-4) & 2.02970 &  9.73319(-5)& 2.02090 & 3.63870(-3) & 0.94265  \\ \hline
\end{array}
$$
\caption{Example 3: Errors and rate of convergences in $\|u-u_h\|, \; \|u - u_h\|_{\infty}$ and $\|u-u_h\|_1$.}
\label{ex-t3}
\end{table}

\begin{figure}[htb!]
  \centering
  \includegraphics[width=.5\linewidth]{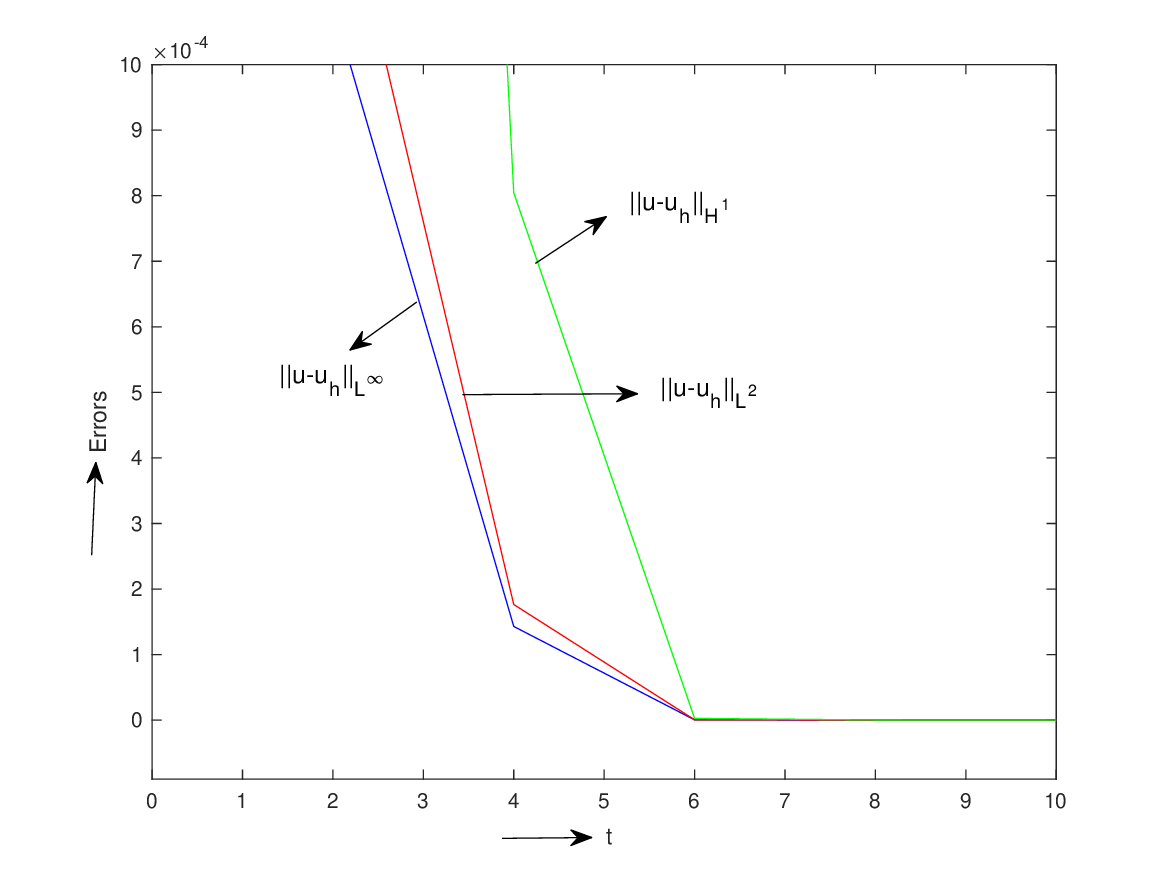}
  \caption{Example 4: The decay estimates in $L^{\infty}, \, L^2$ and $H^{1}$-norms for $\alpha = 1$.}
  \label{ex-f4}
\end{figure}

\noindent
{\bf Example 4.} For the  semilinear weakly damped wave equation, see  \cite{Jos} 
\beas
u_{tt} + \alpha \, u_{t} - \Delta u + f(u) = 0, \; (x_1, x_2) \in \Omega = (0,1) \times (0, 1), \; t > 0
\eeas
with initial conditions
\beas
u(x_1, x_2, 0) = \sin (\pi x_1) \sin (\pi x_2), \quad u_t(x_1, x_2, 0) = -\pi \sin (\pi x_1) \sin (\pi x_2), \; (x_1, x_2) \in \Omega
\eeas
and homogeneous Dirichlet  boundary condition,
where $f(u) = u^3 - u$, we observe that errors decay exponentially in Figure \ref{ex-f4}. 

\begin{figure}[htb!]
\centering
\begin{subfigure}{.5\textwidth}
  \centering
  \includegraphics[width=.95\linewidth]{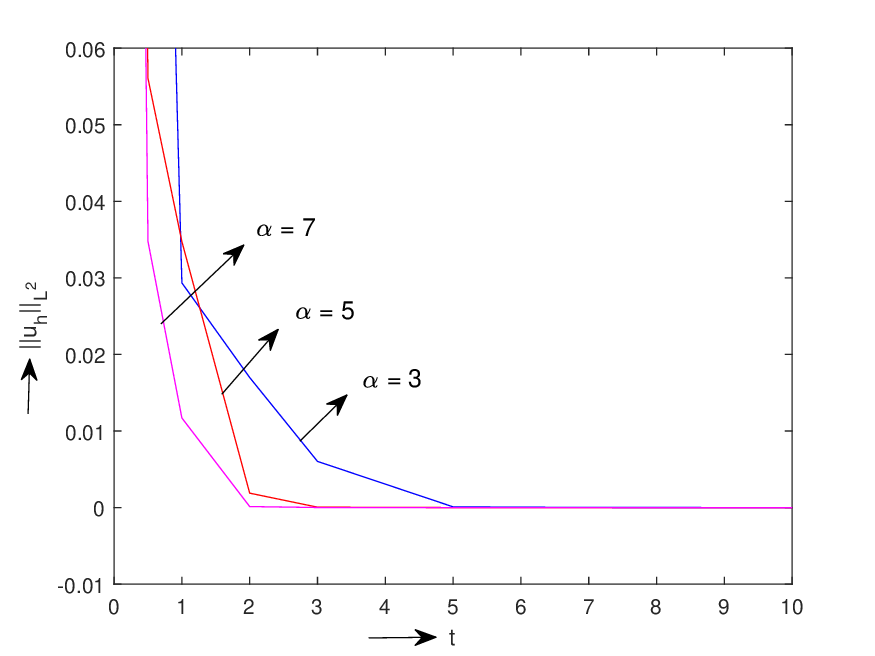}
  \caption{The decay estimates in $L^{2}$-norm.}
  \label{fig:Ex6.1}
\end{subfigure}%
\begin{subfigure}{.5\textwidth}
  \centering
  \includegraphics[width=.95\linewidth]{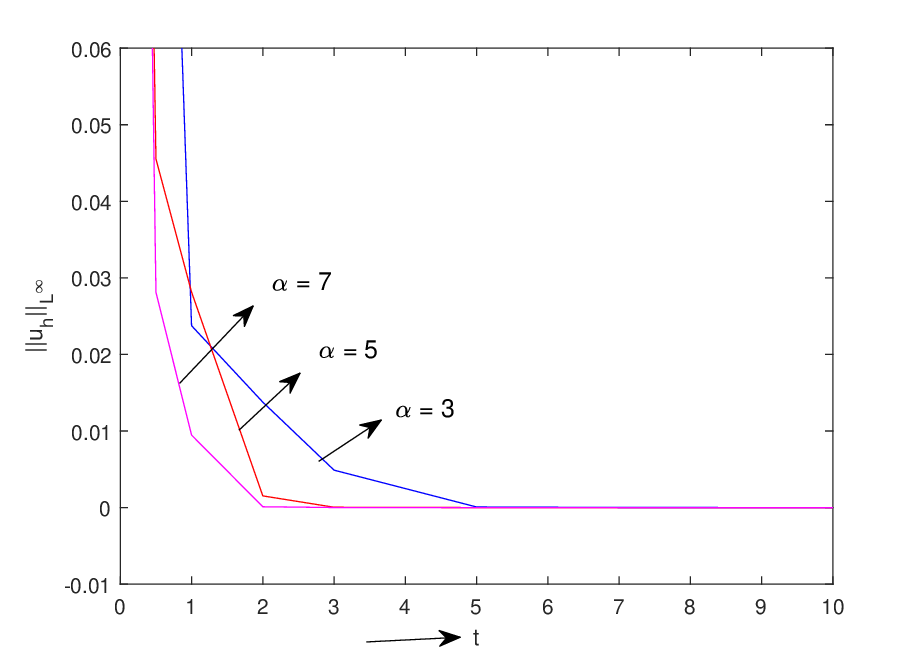}
  \caption{The decay estimates in $L^{\infty}$-norm.}
  \label{fig:Ex6.2}
\end{subfigure}
\caption{Example 4: The decay estimates in $L^2$ and $L^{\infty}$-norms for different $\alpha$ values.}
\label{ex-f4-2}
\end{figure}

The Figure \ref{ex-f4-2} shows that  the decay plots for different values of damping coefficient $\alpha$. It is  observed that $\alpha = 7$ decay exponentially faster than $\alpha = 3$ and $\alpha = 5$. This confirms that exponentially decay phenomenon for the norms $L^2$ and $L^{\infty}$. \\

\begin{figure}[htb!]
\centering
\begin{subfigure}{.45\textwidth}
  \centering
  \includegraphics[width=.95\linewidth]{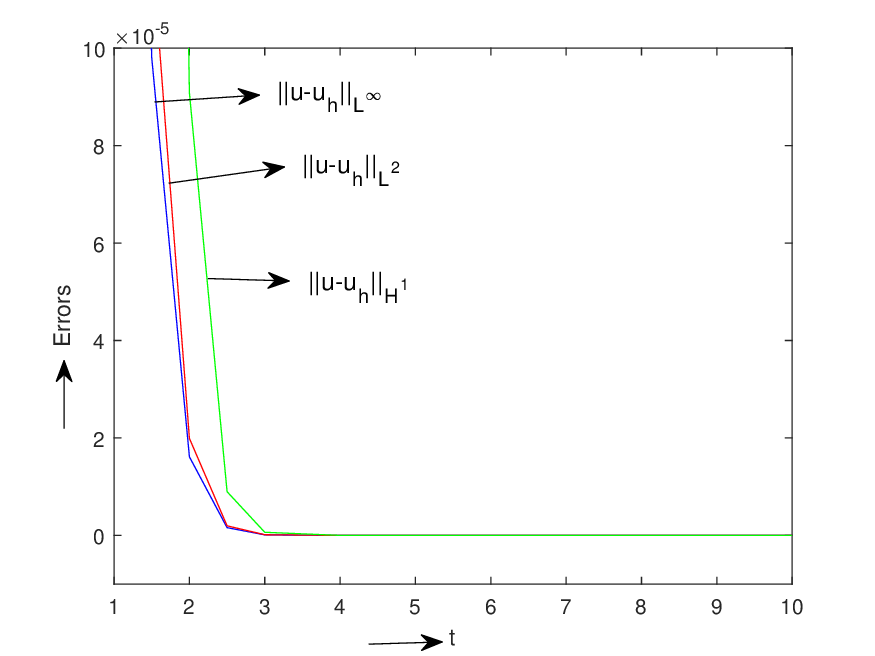}
  \caption{The decay estimate for $\alpha = 10$ and $\beta = 32$.}
  \label{fig:Ex8.1}
\end{subfigure}%
\begin{subfigure}{.45\textwidth}
  \centering
  \includegraphics[width=.95\linewidth]{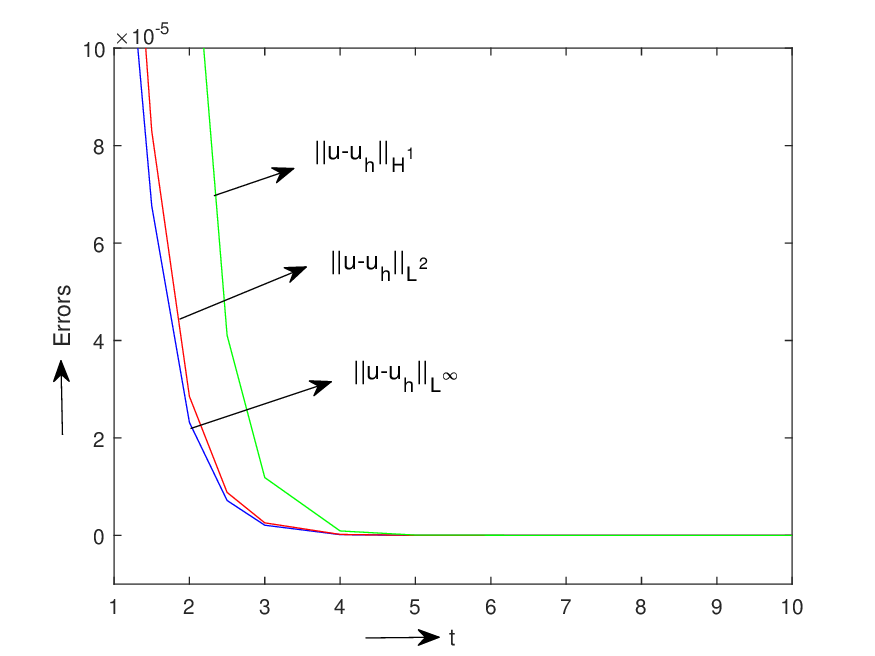}
  \caption{The decay estimate for $\alpha = 10$ and $\beta = 0$.}
  \label{fig:Ex8.2}
\end{subfigure}
\centering
\begin{subfigure}{.45\textwidth}
  \centering
  \includegraphics[width=.95\linewidth]{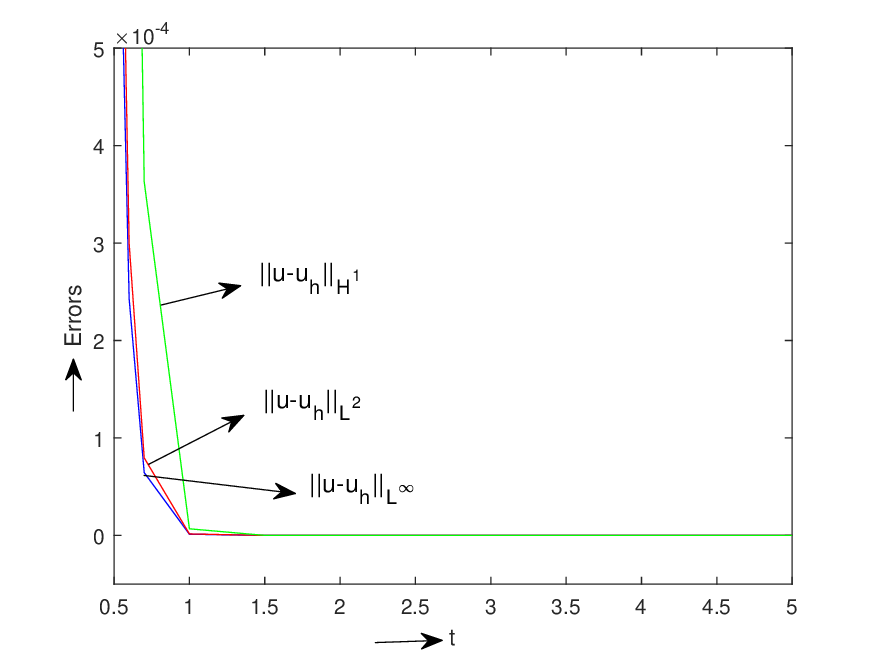}
  \caption{The decay estimate for $\alpha = 40$ and $\beta = 335$.}
  \label{fig:Ex8.3}
\end{subfigure}%
\begin{subfigure}{.45\textwidth}
  \centering
  \includegraphics[width=.95\linewidth]{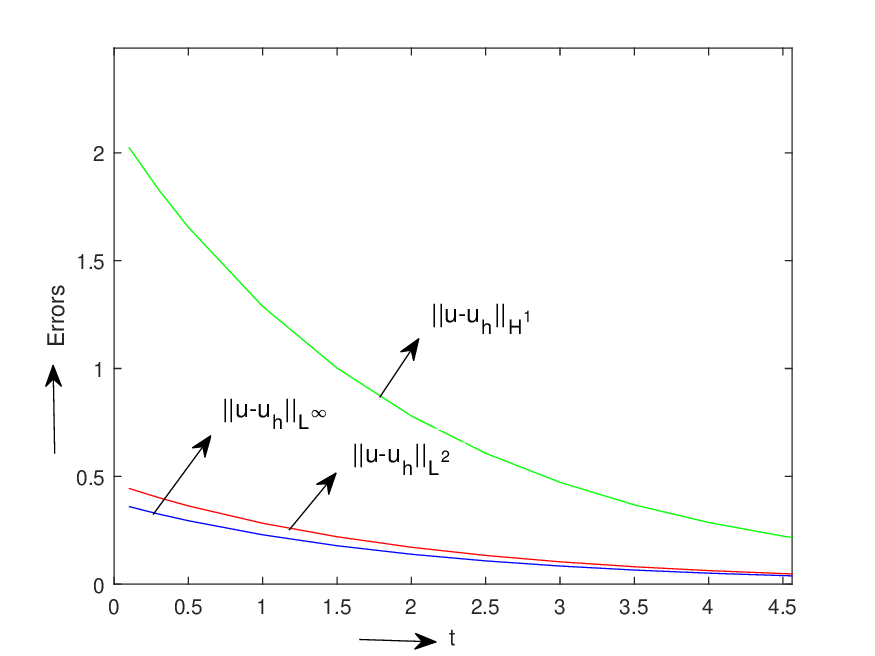}
  \caption{The decay estimate for $\alpha = 40$ and $\beta = 0$.}
  \label{fig:Ex8.4}
\end{subfigure}
\caption{Example 5: The decay estimates in $L^{\infty}, \, L^2$ and $H^{1}$-norms for different $\alpha$ and $\beta$ values.}
\label{ex-f5-1}
\end{figure}

\noindent
{\bf Example 5.}  For the wave equation with viscous damping and compensation
\beas
u_{tt}+ \alpha \, u_{t} + \beta \, u - \Delta u = 0, \; (x_1, x_2) \in \Omega = (0,1) \times (0, 1), \; t > 0
\eeas
with initial conditions
\beas
u(x_1, x_2, 0) = \sin (\pi x_1) \sin (\pi x_2), \quad u_{t}(x_1, x_2, 0) = -\pi \sin (\pi x_1) \sin (\pi x_2), \; (x_1, x_2) \in \Omega
\eeas
and homogeneous  boundary condition,
we calculate the values of damping coefficient $\alpha$ and compensation coefficient $\beta$ from (\ref{condition:delta}). If we choose $\delta = 2$ and $\delta = 5$, that is, decay rate $1$ and $5/2,$ respectively, then we obtain $\alpha = 10, \; \beta = 32$ and $\alpha = 40, \; \beta = 335$, respectively. 
Now the  decay plots for different values of damping coefficient $\alpha$ and compensation coefficient $\beta$ are shown in Figure \ref{ex-f5-1}.\\

\begin{figure}[htb!]
\centering
\begin{subfigure}{.45\textwidth}
  \centering
  \includegraphics[width=.9\linewidth]{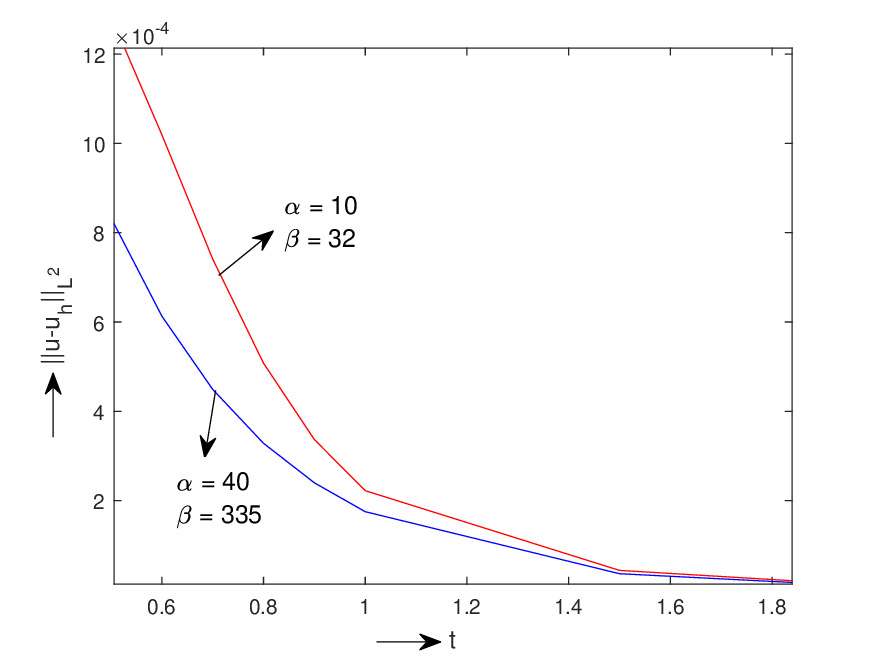}
  \caption{$\alpha = 10, 40$ and $\beta = 32, 335$ in $L^2$-norm.}
  \label{fig:Ex9.3}
\end{subfigure}%
\begin{subfigure}{.45\textwidth}
  \centering
  \includegraphics[width=.9\linewidth]{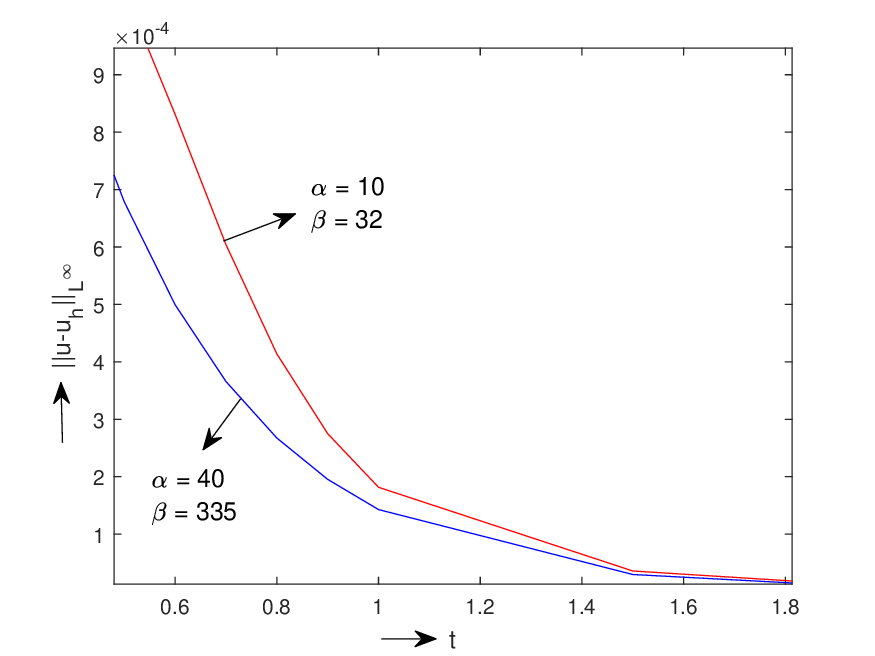}
  \caption{$\alpha = 10, 40$ and $\beta = 32, 335$  in $L^{\infty}$-norm}
  \label{fig:Ex9.4}
\end{subfigure}
\caption{Example 5: The decay estimates for different pairs of damping and compensation coefficients.}
\label{ex-f5-2}
\end{figure}

In Figure \ref{ex-f5-2}, we compute the decay rates for different values of damping and compensation parameters.
\noindent 

\begin{figure}[htb!]
\centering
\begin{subfigure}{.45\textwidth}
  \centering
  \includegraphics[width=.95\linewidth]{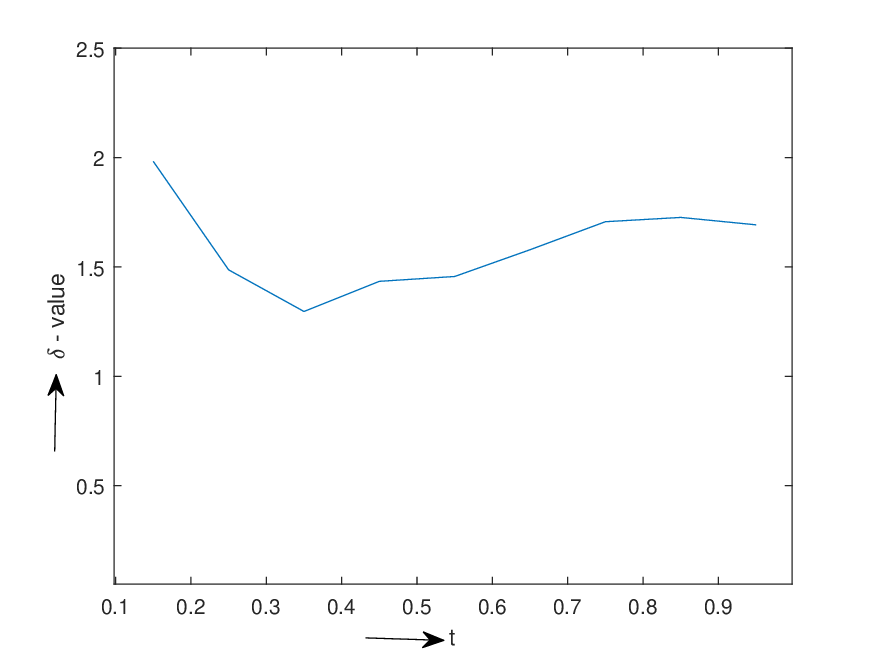}
  \caption{Computation of $\delta$ when $\alpha = 10$ and $\beta = 32$.}
  \label{fig:Ex6.5}
\end{subfigure}%
\begin{subfigure}{.45\textwidth}
  \centering
  \includegraphics[width=.95\linewidth]{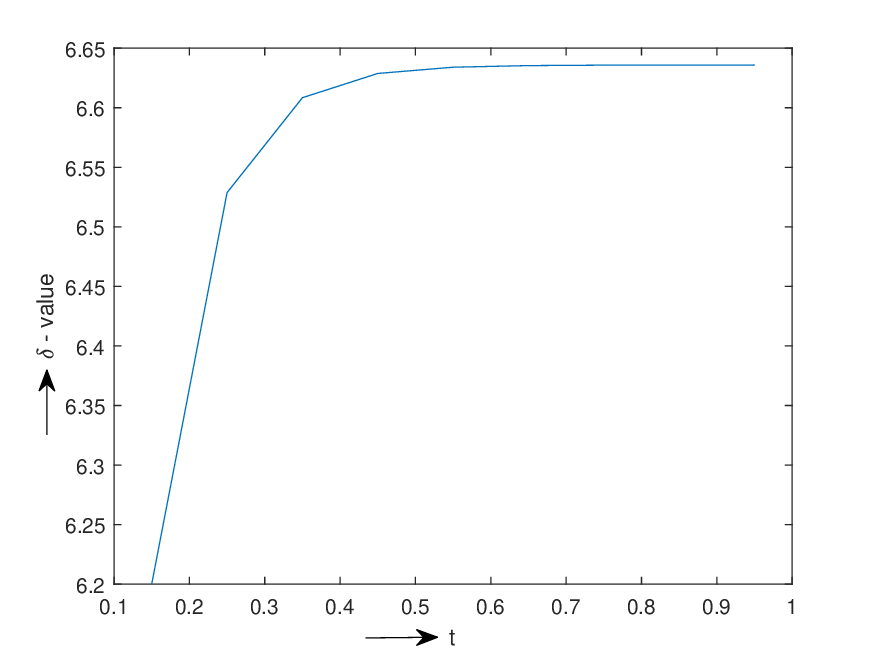}
  \caption{Computation of $\delta$ when $\alpha = 40$ and $\beta = 335$.}
  \label{fig:Ex6.6}
\end{subfigure}
\caption{Example 5: Computation of $\delta$ numerically.}
\label{ex-f5-3}
\end{figure}

Below, in Figure \ref{ex-f5-3}, we compute the decay rate numerically for different values of damping and compensation parameters.

\noindent

{\bf{Observations}}:
\begin{itemize}
\item From Figure \ref{ex-f5-1}(a), (c), it is noted that when $\alpha = 40$ with $\beta = 335$, the errors decay exponentially  with rate $\delta/2=5/2$ faster than $\alpha = 10$ and $\beta = 32$ with rate $\delta/2=1.$ 
This confirms that for any arbitrary  $\delta,$  one may choose damping coefficient $\alpha$ and compensation parameter $\beta$ appropriately so that  errors in $L^2, H^1$ and $L^\infty$-norms  decay exponentially with decay rate $\delta/2.$ 
Now, we examine the decay estimates by setting compensation coefficient $\beta = 0$. 
Comparing both the decay estimates, it is observed  that the errors in Figure \ref{ex-f5-1}(a),(c) decay exponentially faster than the errors in Figure \ref{ex-f5-1}(b),(d). This confirms that the compensation term $\beta$ is helping in the weakly damping equation to get the errors decay exponentially faster.

\item It is further observe through  numerical experiments for the wave equation with different viscous damping coefficients and compensation coefficients in Figure \ref{ex-f5-2} that for large decay rates one may choose the  compensation term and damping coefficient large as given in the subsection 4.3 which is better than  the decay rate than decay predicted in the Sections 2 and 3. Say, for example with decay rates $1$,  the damping coefficients $\alpha=10$ and the compensation parameter $32,$   the predicted decay rate  as in Sections 2 and 3 for the Example 5 is less than equal to 
$\frac{1}{3} \min( 10, 33/20)=11/20,$  with $\lambda_1=1$ which confirms our results in subsection 4.3.

\item Figure \ref{ex-f5-3} shows the calculation of the decay rate $\delta/2$ numerically. When $\alpha = 10, \; \beta = 32$,  it is noticed that $\delta$ is converging close to $2,$ that is, decay rate $1,$ which confirms theoretical result in subsection 4.3. Further with  $\alpha = 40, \; \beta = 335,$ $\delta$ is converging close to $6.6,$  that is, the decay rate in this case is roughly $3.2$, which seems to be better than the decay rate $2.5$ as predicted by the Theorem 4.2. This suggests that the choice of $\alpha$ and $\beta$ in terms of $\delta$ may not be conservative. 
\end{itemize}

\subsection{Conclusions}
In this article, the uniform exponential decay estimates for the linear weakly damped wave equation are developed and analyzed  for continuous and semidiscrete problem. Semidiscrete approximations are obtained by applying FEM to discretize  in space directions keeping the time variable continuous. Compared to the existing literature, improved decay rates with rates lying in a range are derived. 
It is further observed that optimal error estimates, which  depict the decay behaviour are proved with minimal smoothness assumptions on the initial data.The present analysis  is extended  to problems with inhomogeneous forcing function, space dependent damping coefficient, viscous damping and compensation. As a consequence of our abstract analysis, the proof technique is also  generalized to a  weakly damped beam equation. Several numerical experiments are performed to validate the theoretical results established in this article. The optimal rate of convergence is achieved in Table 1-5 and uniform exponential decay behaviour is observed in Figure 1-8. 
In examples 5-6, it  is  shown numerically that the semidiscrete solution of the semilinear weakly damped equation decays exponentially, and  in future, we shall develop similar results as in linear case. Moreover,  our future investigation will include  the uniform exponential decay estimates for the complete discrete schemes. 

~\\
{\bf Funding:} The authors would like to thank SERB, India for the financial support received (Grant Number CRG/2023/000721).

\addcontentsline{toc}{section}{References}
\bibliographystyle{plain}
\bibliography{reference.bib}

\begin{thebibliography}{10}

\bibitem{B-1976}
Garth~A Baker.
\newblock Error estimates for finite element methods for second order
  hyperbolic equations.
\newblock {\em SIAM Journal on Numerical Analysis}, 13(4):564--576, 1976.

\bibitem{JMB}
John~M Ball.
\newblock Global attractors for damped semilinear wave equations.
\newblock {\em Discrete and Continuous Dynamical Systems}, 10(1/2):31--52,
  2004.

\bibitem{Bre}
Susanne~C Brenner and L~Ridgway Scott.
\newblock {\em The mathematical theory of finite element methods}, volume~15.
\newblock 2008.

\bibitem{chen-2}
G~Chen, SA~Fulling, FJ~Narcowich, and C~Qi.
\newblock An asymptotic average decay rate for the wave equation with variable
  coefficient viscous damping.
\newblock {\em SIAM Journal on Applied Mathematics}, 50(5):1341--1347, 1990.

\bibitem{chen-1}
Goong Chen.
\newblock Control and stabilization for the wave equation in a bounded domain.
\newblock {\em SIAM Journal on Control and Optimization}, 17(1):66--81, 1979.

\bibitem{chen-3}
Goong Chen.
\newblock Control and stabilization for the wave equation in a bounded domain,
  part ii.
\newblock {\em SIAM Journal on Control and Optimization}, 19(1):114--122, 1981.

\bibitem{CZ}
Steven Cox and Enrique Zuazua.
\newblock The rate at which energy decays in a damped string.
\newblock {\em Communications in Partial Differential Equations},
  19(1-2):213--243, 1994.

\bibitem{HET}
Herbert Egger and Thomas Kugler.
\newblock {\em Uniform Exponential Stability of Galerkin Approximations for a
  Damped Wave System}, pages 107--129.
\newblock Springer International Publishing, Cham, 2019.

\bibitem{LCE}
Lawrence~C Evans.
\newblock {\em Partial differential equations}, volume~19.
\newblock American Mathematical Society, 2022.

\bibitem{GR}
Jerome~A Goldstein and Steven~I Rosencrans.
\newblock Energy decay and partition for dissipative wave equations.
\newblock {\em Journal of Differential Equations}, 36(1):66--73, 1980.

\bibitem{ff++}
Fr{\'e}d{\'e}ric Hecht.
\newblock New development in freefem++.
\newblock {\em Journal of Numerical Mathematics}, 20(3-4):251--266, 2012.

\bibitem{Jos}
Joseph~W Jerome.
\newblock The multidimensional damped wave equation: maximal weak solutions for
  nonlinear forcing via semigroups and approximation.
\newblock {\em Numerical Functional Analysis and Optimization},
  41(16):1970--1989, 2020.

\bibitem{Kap}
Lev Kapitanski.
\newblock Minimal compact global attractor for a damped semilinear wave
  equation.
\newblock {\em Communications in Partial Differential Equations},
  20(7-8):1303--1323, 1995.

\bibitem{KPY}
Samir Karaa, Amiya~K Pani, and Sangita Yadav.
\newblock A priori hp-estimates for discontinuous galerkin approximations to
  linear hyperbolic integro-differential equations.
\newblock {\em Applied Numerical Mathematics}, 96:1--23, 2015.

\bibitem{Lions-Strauss}
Jacques-Louis Lions and Walter~A Strauss.
\newblock Some non-linear evolution equations.
\newblock {\em Bulletin de la Soci{\'e}t{\'e} Math{\'e}matique de France},
  93:43--96, 1965.

\bibitem{RTT}
Karim Ramdani, Tak{\'e}o Takahashi, and Marius Tucsnak.
\newblock Uniformly exponentially stable approximations for a class of second
  order evolution equations: Application to lqr problems.
\newblock {\em ESAIM: Control, Optimisation and Calculus of Variations},
  13(3):503--527, 2007.

\bibitem{Rau}
Jeffrey Rauch.
\newblock Qualitative behavior of dissipative wave equations on bounded
  domains.
\newblock {\em Archive for Rational Mechanics and Analysis}, 62(1):77--85,
  1976.

\bibitem{Rauch}
Jeffrey Rauch.
\newblock On convergence of the finite element method for the wave equation.
\newblock {\em SIAM Journal on Numerical Analysis}, 22(2):245--249, 1985.

\bibitem{Russell-1}
David~L Russell.
\newblock Decay rates for weakly damped systems in hilbert space obtained with
  control-theoretic methods.
\newblock {\em Journal of Differential Equations}, 19(2):344--370, 1975.

\bibitem{S}
Rajen~K Sinha.
\newblock Finite element approximations with quadrature for second-order
  hyperbolic equations.
\newblock {\em Numerical Methods for Partial Differential Equations: An
  International Journal}, 18(4):537--559, 2002.

\bibitem{SP}
Rajen~K Sinha and Amiya~K Pani.
\newblock The effect of spatial quadrature on finite element galerkin
  approximations to hyperbolic integro-differential equations.
\newblock {\em Numerical Functional Analysis and Optimization},
  19(9-10):1129--1153, 1998.

\bibitem{Temam}
Roger Temam.
\newblock {\em Infinite-dimensional dynamical systems in mechanics and
  physics}, volume~68.
\newblock Springer Science \& Business Media, 2012.

\bibitem{VT}
Vidar Thom{\'e}e.
\newblock {\em Galerkin finite element methods for parabolic problems},
  volume~25.
\newblock Springer Science \& Business Media, 2007.

\bibitem{BV}
{Babin, Anatoli{\u\i}} Vladimirovich and Mark~I Vishik.
\newblock {\em Attractors of evolution equations}, volume~25.
\newblock Elsevier, 1992.

\bibitem{Z-2005}
Enrique Zuazua.
\newblock Propagation, observation, and control of waves approximated by finite
  difference methods.
\newblock {\em SIAM Review}, 47(2):197--243, 2005.

\end{thebibliography}
\end{document}